\documentstyle[12pt]{amsart}
\begin{document}
\title{Classifications and isolation phenomena of Bi-Harmonic Maps and Bi-Yang-Mills Fields}
\title[Bi-Harmonic Maps and Bi-Yang-Mills Fields]
{Classification and isolation phenomena of Bi-Harmonic Maps and Bi-Yang-Mills Fields}
\author{Toshiyuki  Ichiyama}
\address{Faculty of Economics, Asia University, Sakai 5-24-10, Musashino, Tokyo, 180-8624, Japan}
\email{ichiyama@@asia-u.ac.jp}
\author{Jun-ichi Inoguchi} 
\address{Department of Mathematics, Faculty of Education, Utsunomiya University, Utsunomiya, 321-8505, Japan} 
\email{inoguchi@@cc.utsunomiya-u.ac.jp}
\author{Hajime Urakawa}
\address{Division of Mathematics, Graduate School of Information Sciences, Tohoku University, Aoba 6-3-09, Sendai, 980-8579, Japan}
\email{urakawa@@math.is.tohoku.ac.jp}
\maketitle
\begin{abstract}
Classifications of all biharmonic isoparametric hypersurfaces in the unit sphere, and all biharmonic homogeneous real hypersurfaces in the complex or quaternionic projective spaces are shown. Answers in case of bounded geometry to Chen's conjecture or Caddeo, Montaldo and Piu's one on biharmonic maps into a space of  non positive curvature are given. Gauge field analogue is shown, indeed, 
the isolation phenomena of bi-Yang-Mills fields are obtained. 
\end{abstract}
\numberwithin{equation}{section}
\theoremstyle{plain}
\newtheorem{df}{Definition}[section]
\newtheorem{th}{Theorem}[section]
\newtheorem{prop}{Proposition}[section]
\newtheorem{lem}{Lemma}[section]
\newtheorem{cor}{Corollary}[section]
\section{Introduction}
Theory of harmonic maps plays a central roll in variational problems, which are by definition for smooth maps 
between Riemannian manifolds $\varphi:\,M\rightarrow N$, 
critical maps of the energy functional 
$E(\varphi)=\frac12\int_M\Vert d\varphi\Vert^2\,v_g$. 
By extending the notion of harmonic maps, in 1983, J. Eells and L. Lemaire \cite{EL1} proposed the problem 
to consider the $k$-harmonic maps which are critical maps of the functional
$$
E_k(\varphi)=\frac12\int_M
\Vert(d+\delta)^k\varphi\Vert^2\,v_g,\quad (k=1,2,\cdots).
$$
After G.Y. Jiang \cite{J} studied the first and second variation formulas of $E_k$ for $k=2$, whose critical maps are called biharmonic maps, there have been extensive studies in this area (for instance, see \cite{CMP}, \cite{LO}, \cite{LO2}, \cite{O1}, \cite{MO1}, 
\cite{I}, \cite{II}, \cite{S1}, etc.). Harmonic maps are always biharmonic maps by definition. One of main central problems is to classify the biharmonic maps, or to ask whether or not the converse to the above is true when the target Riemannian manifold $(N,h)$ is non positive curvature (B. Y. Chen's conjecture \cite{C} or Caddeo, Montaldo and Piu's one \cite{CMP}). 
In this paper, (1) we classify all biharmonic hypersurfaces isoparametric hypersurfaces in the unit sphere, i.e., whose principal curvatures are constant, in \S 3, 4, and (2) we give the first examples and classify all biharmonic homogeneous real hypersurfaces in the complex or quaternionic projective spaces in \S 5, 6, 7. 
Next, we give answers to Chen's conjecture and Caddeo, Montaldo and Piu's one in \S 8. Indeed, we show all biharmonic maps or biharmonic submanifolds of bounded geometry into the target space which is  non positive curvature, 
must be harmonic. 
Here, that biharmonic maps are of bounded geometry means that the curvature of the domain manifold is bounded, and the norms of  the tension field and its covariant derivative are $L^2$. 
\par
Recently, the notion of gauge field analogue of biharmonic maps, i.e., bi-Yang-Mills fields was proposed (\cite{BU}). 
In this paper, we show the isolation phenomena of bi-Yang-Mills fields like the one for Yang-Mills fields (cf. Bourguignon-Lawson \cite{BL}), i.e., all bi-Yang-Mills fields over compact Riemanian manifolds of which Ricci curvature are bounded below by a positive constant $k$,  and the pointwise norm of curvature tensor are bounded above by $k/2$, must be Yang-Mills fields. We also show the $L^2$-isolation phenopmena which are similar as 
Min-Oo's result (\cite{MO}) for Yang-Mills fields.  
 These interesting phenomena can be regarded that the similar phenomena as the biharmonic maps should occur.
\section{Preliminaries}
In this section, we prepare materials for the first variation formula for the bi-energy functional and bi-harmonic maps. 
Let us recall the definition of a harmonic map $\varphi:\,(M,g)\rightarrow (N,h)$, of a comoact Riemannian manifold $(M,g)$ into another Riemannian manifold $(N,h)$, 
which is an extremal 
of the {\em energy functional} defined by 
$$
E(\varphi)=\int_Me(\varphi)\,v_g, 
$$
where $e(\varphi):=\frac12\vert d\varphi\vert^2$ is called the energy density 
of $\varphi$.  
That is, for all variation $\{\varphi_t\}$ of $\varphi$ with 
$\varphi_0=\varphi$, 
\begin{equation}
\frac{d}{dt}\bigg\vert_{t=0}E(\varphi_t)=-\int_Mh(\tau(\varphi),V)v_g=0,
\end{equation}
where $V\in \Gamma(\varphi^{-1}TN)$ is a variation vector field along $\varphi$ which is given by 
$V(x)=\frac{d}{dt}\vert_{t=0}\varphi_t(x)\in T_{\varphi(x)}N$
$(x\in M)$, 
and the {\em tension field} of $\varphi$ is given by 
$\tau(\varphi)
=\sum_{i=1}^mB(\varphi)(e_i,e_i)\in \Gamma(\varphi^{-1}TN)$,  
where 
$\{e_i\}_{i=1}^m$ is a locally defined frame field on $(M,g)$.   
The second fundamental form  $B(\varphi)$ of $\varphi$ is  
defined by 
\begin{align}
B(\varphi)(X,Y)&=(\widetilde{\nabla} d \varphi)(X,Y)\nonumber\\
&=(\widetilde{\nabla}_Xd\varphi)(Y)\nonumber\\
&=\overline{\nabla}_X(d\varphi(Y))-d\varphi(\nabla_XY)\nonumber\\
&=\nabla^N_{d\varphi(X)}d\varphi(Y)-d\varphi(\nabla_XY),
\end{align}
for all vector fields $X, Y\in {\frak X}(M)$. 
Furthermore, 
$\nabla$, and
$\nabla^N$, 
 are connections on $TM$, $TN$  of $(M,g)$, $(N,h)$, respectively, and 
$\overline{\nabla}$, and $\widetilde{\nabla}$ are the induced one on $\varphi^{-1}TN$, and $T^{\ast}M\otimes \varphi^{-1}TN$, respectively. By (2.1), $\varphi$ is harmonic if and only if $\tau(\varphi)=0$. 
\par
The second variation formula of the energy functional is also well known which is given as follows. Assume that 
$\varphi$ is harmonic. 
Then, 
\begin{equation}
\frac{d^2}{dt^2}\bigg\vert_{t=0}E(\varphi_t)
=\int_Mh(J(V),V)v_g, 
\end{equation}
where 
$J$ is an ellptic differential operator, called 
{\em Jacobi operator}  acting on 
$\Gamma(\varphi^{-1}TN)$ given by 
\begin{equation}
J(V)=\overline{\Delta}V-{\mathcal R}(V),
\end{equation}
where 
$\overline{\Delta}V=\overline{\nabla}^{\ast}\overline{\nabla}V$ 
is the {\em rough Laplacian} and 
${\mathcal R}$ is a linear operator on $\Gamma(\varphi^{-1}TN)$
given by 
${\mathcal R}V=
\sum_{i=1}^mR^N(V,d\varphi(e_i))d\varphi(e_i)$,
and $R^N$ is the curvature tensor of $(N,h)$ given by 
$R^N(U,V)=\nabla^N{}_U\nabla^N{}_V-\nabla^N{}_V\nabla^N{}_U-\nabla^N{}_{[U,V]}$ for $U,\,V\in {\frak X}(N)$.   
\par
J. Eells and L. Lemaire proposed (\cite{EL1}) polyharmonic ($k$-harmonic) maps and 
Jiang studied (\cite{J}) 
the first and second variation formulas of bi-harmonic maps. Let us consider the {\em bi-energy functional} 
defined by 
\begin{equation}
E_2(\varphi)=\frac12\int_M\vert\tau(\varphi)\vert ^2v_g, 
\end{equation}
where 
$\vert V\vert^2=h(V,V)$, $V\in \Gamma(\varphi^{-1}TN)$.  
Then, the first variation formula is given as follows. 
\begin{th}
\quad $($the first variation formula$)$ 
\begin{equation}
\frac{d}{dt}\bigg\vert_{t=0}E_2(\varphi_t)
=-\int_Mh(\tau_2(\varphi),V)v_g,
\end{equation}
where 
\begin{equation}
\tau_2(\varphi)=J(\tau(\varphi))=\overline{\Delta}\tau(\varphi)-{\mathcal R}(\tau(\varphi)),
\end{equation}
$J$ is given in $(2.4)$. 
\end{th}
\vskip0.6cm\par
For the second variational formula, see \cite{J} or \cite{IIU}.
\vskip0.6cm\par
\begin{df}
A smooth map $\varphi$ of $M$ into $N$ is called to be 
{\em bi-harmonic} if 
$\tau_2(\varphi)=0$. 
\end{df} 
\vskip0.6cm\par
For later use, we need the following three lemmas. 
\begin{lem} $($Jiang$)$ 
Let $\varphi:\,(M^m,g)\rightarrow (N^n,h)$ be an isometric immersion of which mean curvuture vector field ${\Bbb H}=\frac{1}{m}\tau(\varphi)$ is parallel, i.e., $\nabla^{\perp}{\Bbb H}=0$, 
where $\nabla^{\perp}$ is the induced connection of the normal bundle $T^{\perp}M$ by $\varphi$. Then, 
\begin{align}
\overline{\Delta}\tau(\varphi)
&=
\sum_{i=1}^mh(\overline{\Delta}\tau(\varphi),d\varphi(e_i))d\varphi(e_i)\nonumber\\
&\quad
-\sum_{i,j=1}^mh(\overline{\nabla}_{e_i}\tau(\varphi),d\varphi(e_j))
(\widetilde{\nabla}_{e_i}d\varphi)(e_j),
\end{align}
where $\{e_i\}$ is a locally defined orthonormal frame field 
of $(M,g)$. 
\end{lem} 
\begin{pf}
Let us recall the definition of $\nabla^{\perp}$: 
For any section $\xi\in \Gamma(T^{\perp}M)$, we decompose $\overline{\nabla}_X\xi$ 
according to $TN\vert_M=TM\oplus T^{\perp}M$ as follows. 
$$
\overline{\nabla}_X\xi=\nabla^N_{\varphi_{\ast}X}\xi
=\nabla^{T}_{\varphi_{\ast}X}\xi+\nabla^{\perp}_{\varphi_{\ast}X}\xi.
$$
By the assumption $\nabla^{\perp}{\Bbb H}=0$, 
i.e., $\nabla_{\varphi_{\ast}X}^{\perp}\tau(\varphi)=0$ for all 
$X\in {\frak X}(M)$, 
we have 
\begin{equation}
\overline{\nabla}_X\tau(\varphi)=
\nabla^T_{\varphi_{\ast}X}\tau(\varphi)\in \Gamma(\varphi_{\ast}TM).
\end{equation}
Thus, for all $i=1,\cdots,m$, 
\begin{equation}
\overline{\nabla}_{e_i}\tau(\varphi)=
\sum_{j=1}^mh(\overline{\nabla}_{e_i}\tau(\varphi),d\varphi(e_j))
d\varphi(e_j)
\end{equation}
because $\{d\varphi(e_j)_x\}_{j=1}^m$ is an orthonormal basis  with respect to $h$, of 
$\varphi_{\ast}T_xM$ $(x\in M)$. 
\par
Now let us calculate 
\begin{equation}
\overline{\nabla}^{\ast}\overline{\nabla}\tau(\varphi)
=-\sum_{i=1}^m
\{
\overline{\nabla}_{e_i}\overline{\nabla}_{e_i}\tau(\varphi)
-\overline{\nabla}_{\nabla_{e_i}e_i}\tau(\varphi)
\}.
\end{equation}
Indeed, we have 
\begin{align}
\overline{\nabla}_{e_i}\overline{\nabla}_{e_i}\tau(\varphi)
&=\sum_{j=1}^m\{
h(\overline{\nabla}_{e_i}\overline{\nabla}_{e_i}\tau(\varphi))
+h(\overline{\nabla}_{e_i}\tau(\varphi),\overline{\nabla}_{e_i}
d\varphi(e_j))
\}d\varphi(e_j)\nonumber\\
&+\sum_{j=1}^m
h(\overline{\nabla}_{e_i}\tau(\varphi),d\varphi(e_j))\overline{\nabla}_{e_i}d\varphi(e_i), 
\end{align}
and 
\begin{equation}
\overline{\nabla}_{\nabla_{e_i}e_i}\tau(\varphi)=
\sum_{j=1}^m
h(\overline{\nabla}_{{\nabla}_{e_i}e_i},d\varphi(e_j))d\varphi(e_j),
\end{equation}
so that we have 
\begin{align}
\overline{\nabla}^{\ast}\overline{\nabla}
\tau(\varphi)&=
\sum_{j=1}^mh(\overline{\nabla}^{\ast}\overline{\nabla}\tau(\varphi),
d\varphi(e_j))d\varphi(e_j)\nonumber\\
&-\sum_{i,j=1}^m
\{
h(\overline{\nabla}_{e_i}\tau(\varphi),\overline{\nabla}_{e_i}
d\varphi(e_j))
\}d\varphi(e_j)\nonumber\\
&\qquad\quad+
h(\overline{\nabla}_{e_i}\tau(\varphi),d\varphi(e_j))\overline{\nabla}_{e_i}d\varphi(e_i)
\}.
\end{align}
Denoting $\nabla_{e_i}e_j=\sum_{k=1}^m\Gamma_{ij}^ke_k$, we have $\Gamma^k_{ij}+\Gamma^j_{ik}=0$. 
Since 
$(\widetilde{\nabla}_{e_i}d\varphi)(e_j)\newline=
\overline{\nabla}_{e_i}(d\varphi(e_j))-d\varphi(\nabla_{e_i}e_j)$ 
is a local section of  $T^{\perp}M$, 
we have for the the second term of the RHS of (2.14), 
for each fixed 
$i=1,\cdots,m$, 
\begin{align}
\sum_{j=1}^m
h(\overline{\nabla}_{e_i}\tau(\varphi),
&\overline{\nabla}_{e_i}d\varphi(e_j))d\varphi(e_j)\nonumber\\
&=
\sum_{j=1}^m
h(\overline{\nabla}_{e_i}\tau(\varphi),
(\widetilde{\nabla}_{e_i}d\varphi)(e_j)
+d\varphi(\nabla_{e_i}e_j))d\varphi(e_j)
\nonumber\\
&=
\sum_{j=1}^m
h(\overline{\nabla}_{e_i}\tau(\varphi),
d\varphi(\nabla_{e_i}e_j))d\varphi(e_j)
\nonumber\\
&=\sum_{j,k=1}^m
h(\overline{\nabla}_{e_i}\tau(\varphi),
d\varphi(e_k))\,d\varphi(\Gamma^k_{ij}e_j)
\nonumber\\
&=-\sum_{j,k=1}^m
h(\overline{\nabla}_{e_i}\tau(\varphi),
d\varphi(e_k))\,d\varphi(\Gamma^j_{ik}e_j)
\nonumber\\
&
=-\sum_{k=1}^m
h(\overline{\nabla}_{e_i}\tau(\varphi),
d\varphi(e_k))\,d\varphi(\nabla_{e_i}e_k).
\end{align}
Substituting (2.15) into (2.14), we have the desired (2.8). 
\end{pf}
\vskip0.6cm\par
\begin{lem} $($Jiang$)$
Under the same assumption as Lemma 2.1, we have 
\begin{align}
\overline{\Delta}\tau(\varphi)
&=
-\sum_{j,k=1}^mh(\tau(\varphi),R^N(d\varphi(e_j),d\varphi(e_k))d\varphi(e_k))d\varphi(e_j)
\nonumber\\
&\quad+\sum_{i,j=1}^m
h(\tau(\varphi),(\widetilde{\nabla}_{e_i}d\varphi)(e_j))
(\widetilde{\nabla}_{e_i}d\varphi)(e_j).
\end{align}
\end{lem}
\begin{pf}
Since 
$h(\tau(\varphi),d\varphi(e_j))=0$, differentiating it 
by $e_i$, 
we have 
\begin{align}
h(\overline{\nabla}_{e_i}\tau(\varphi),d\varphi(e_j))
&=-h(\tau(\varphi),\overline{\nabla}_{e_i}d\varphi(e_j))\nonumber\\
&=-h(\tau(\varphi),\overline{\nabla}_{e_i}d\varphi(e_j)
-d\varphi(\nabla_{e_i}e_j))\nonumber\\
&=-h(\tau(\varphi),(\widetilde{\nabla}_{e_i}d\varphi)(e_j)).
\end{align}
\par
For the first term of (2.8), we have for each $j=1,\cdots,m$, 
\begin{align}
h(\overline{\Delta}\tau(\varphi),
d\varphi(e_j))
&-2\sum_{i=1}^m
h(\overline{\nabla}_{e_i}\tau(\varphi),\overline{\nabla}_{e_i}d\varphi(e_j))\nonumber\\
&+h(\tau(\varphi),\overline{\Delta}d\varphi(e_j))=0,
\end{align}
which follows by the expression (2.11) of 
$\overline{\Delta}\tau(\varphi)$,
differentiating the first equation of (2.17) by $e_i$,  and
doing $h(\tau(\varphi),d\varphi(e_j))=0$ 
by $\nabla_{e_i}e_i$. 
\par
For the second term of (2.8), we have by (2.9) and (2.17), 
\begin{align}
h(\overline{\nabla}_{e_i}\tau(\varphi),\overline{\nabla}_{e_i}d\varphi(e_j))
&=
h(\overline{\nabla}_{e_i}\tau(\varphi), 
(\widetilde{\nabla}_{e_i}d\varphi)(e_j)+d\varphi(\nabla_{e_i}e_j))
\nonumber\\
&=
h(\overline{\nabla}_{e_i}\tau(\varphi), d\varphi(\nabla_{e_i}e_j))
\nonumber\\
&=
-h(\tau(\varphi),(\widetilde{\nabla}_{e_i}d\varphi)(\nabla_{e_i}e_j)).
\end{align}
\par
For the third term $h(\tau(\varphi),\overline{\Delta}d\varphi(e_j))$ of (2.18), we have 
\begin{align}
h(\tau(\varphi),\overline{\Delta}d\varphi(e_j))
&=
\sum_{k=1}^m
h(\tau(\varphi),R^N(d\varphi(e_j),d\varphi(e_k))d\varphi(e_k))\nonumber\\
&\quad
-2\sum_{k=1}^m
h(\tau(\varphi),(\widetilde{\nabla}_{e_k}d\varphi)(\nabla_{e_k}e_j)).
\end{align}
Because, 
by making use of   
$(\widetilde{\nabla}_Xd\varphi)(Y)=
\overline{\nabla}_{X}(d\varphi(Y))-d\varphi(\nabla_{X}Y)$ 
and 
$h(\tau(\varphi),d\varphi(X))=0$ 
$(X,Y\in {\frak X}(M))$, 
the LHS of (2.20) coincides with 
\begin{align}
h(\tau(\varphi),
&-\sum_{k=1}^m
\{\overline{\nabla}_{e_k}\overline{\nabla}_{e_k}
-\overline{\nabla}_{\nabla_{e_k}e_k}\}d\varphi(e_j))
\nonumber\\
&=
h(\tau(\varphi),
-\sum_{k=1}^m\{
(\widetilde{\nabla}_{e_k}\widetilde{\nabla}_{e_k}d\varphi)(e_j)
+2(\widetilde{\nabla}_{e_k}d\varphi)(\nabla_{e_k}e_j)
\nonumber\\
&\qquad\qquad\qquad\qquad
-(\widetilde{\nabla}_{\nabla_{e_k}e_k}d\varphi)(e_j)
\})\nonumber\\
&=
h(\tau(\varphi),(\widetilde{\nabla}^{\ast}\widetilde{\nabla}d\varphi)(e_j))
-2h(\tau(\varphi),(\widetilde{\nabla}_{e_k}d\varphi)(\nabla_{e_k}e_j))
\nonumber\\
&=h(\tau(\varphi),\Delta d\varphi(e_j)-Sd\varphi(e_j))
-2h(\tau(\varphi),(\widetilde{\nabla}_{e_k}d\varphi)(\nabla_{e_k}e_j)),
\end{align}
where the last equality follows from the Weitzenb\"ock formula for the Laplacian $\Delta=d\delta+\delta d$ acting on 1-forms on $(M,g)$: 
\begin{equation}
\Delta d\varphi=\widetilde{\nabla}^{\ast}\widetilde{\nabla}d\varphi
+Sd\varphi.
\end{equation}
Here, we have
\begin{align}
Sd\varphi(e_j)
&:=\sum_{k=1}^m(\widetilde{R}(e_k,e_j)d\varphi)(e_k)\nonumber\\
&=\sum_{k=1}^m\{
R^N(d\varphi(e_k),d\varphi(e_j))d\varphi(e_k)
-d\varphi(R^M(e_k,e_j)e_k)\},
\end{align}
and 
\begin{equation}
\Delta d\varphi(e_j)
=d\delta d\varphi(e_j)=-d\tau(\varphi)(e_j)
=-\overline{\nabla}_{e_j}\tau(\varphi).
\end{equation}
Substituting these into (2.24), and using 
$h(\tau(\varphi),d\varphi(X))=0$ 
for all $X\in {\frak X}(M)$, (2.24) coincides with 
$$
\sum_{k=1}^m
\{h(\tau(\varphi),R^N(d\varphi(e_j),d\varphi(e_k))d\varphi(e_k))
-2h(\tau(\varphi),
(\widetilde{\nabla}_{e_k}d\varphi)(\nabla_{e_k}e_j))\},
$$ 
which implies (2.20). 
\par
Substituting (2.19) and (2.20) into  (2.18), we have
\begin{align}
h(\overline{\Delta}\tau(\varphi),d\varphi(e_j))
&=-2\sum_{i=1}^mh(\tau(\varphi),(\widetilde{\nabla}_{e_i}d\varphi)(\nabla_{e_i}e_j))
\nonumber\\
&\quad-\sum_{k=1}^mh(\tau(\varphi),R^N(d\varphi(e_j),d\varphi(e_k))d\varphi(e_k))\nonumber\\
&\quad+2\sum_{k=1}^mh(\tau(\varphi),(\widetilde{\nabla}_{e_k}d\varphi)(\nabla_{e_k}e_j))
\nonumber\\
&=\sum_{k=1}^mh(\tau(\varphi),R^N(d\varphi(e_j),d\varphi(e_k))d\varphi(e_k)).
\end{align}
Substituting (2.19) and (2.25) into (2.8), we have (2.16). 
\end{pf}
\vskip0.6cm\par
\begin{lem}
Let $\varphi:\,(M^m,g)\rightarrow (N^{m+1},h)$ be an isometric immersion which is not harmonic. 
Then, the condition that $\Vert\tau(\varphi)\Vert$ is constant is equivalent to the one that 
\begin{equation}
\overline{\nabla}_X\tau(\varphi)\in \Gamma(\varphi_{\ast}TM), \quad\forall \,X\in {\frak X}(M), 
\end{equation}
that is, the mean curvature tensor is parallel with respect to 
$\nabla^{\perp}$. 
\end{lem}
\begin{pf}Assume that $\varphi$ is not harmonic. 
Then, if $\Vert \tau(\varphi)\Vert$ is constant, 
\begin{equation}
Xh(\tau(\varphi),\tau(\varphi))=2h(\overline{\nabla}_X\tau(\varphi),\tau(\varphi))=0
\end{equation}
for all $X\in {\frak X}(M)$, 
so we have 
$\overline{\nabla}_X\tau(\varphi)\in \Gamma(\varphi_{\ast}TM)$ 
because $\dim M=\dim N-1$ and $\tau(\varphi)\not=0$ everywhere  on $M$.  
The converse is true from the above equality (2.27). 
\end{pf}
\vskip0.6cm\par
\vskip0.6cm\par
\section{Biharmonic maps into the unit sphere}
In this section, we give the classification of all the biharmonic isometrically immersed hypersurfaces of  the unit sphere with constant principal curvatures. 
In order to show it, we need the following theorem.  
\begin{th} 
$($cf. Jiang \cite{J}$)$ Let $\varphi: (M^m,g)\rightarrow S^{m+1}\left(
\frac{1}{\sqrt{c}}
\right)$ be an isometric immersion of an $m$-dimensional compact Riemannian manifold $(M^m,g)$ into the $(m+1)$-dimensional sphere with constant sectional curvature $c>0$. Assume that 
the mean curvature of $\varphi$ is nonzero constant. 
 Then, $\varphi$ is biharmonic if and only if 
square of the pointwise norm of $B(\varphi)$ is constant and $\Vert B(\varphi)\Vert^2 =cm$. 
\end{th}
\begin{pf}
 \quad For completeness, we give a brief proof, here. 
By Lemma 2.3, the condition (2.9) holds under the condition that 
the mean curvature of $\varphi$ is constant. 
So, we may apply Lemmas 2.1 and 2.2. 
\par
Since the curvature tensor $R^N$ of $S^{m+1}\left(\frac{1}{\sqrt{c}}\right)$ is given by 
$$
R^N(U,V)W=c\{h(V,W)U-h(W,U)V\},\quad U,V,W\in {\frak X}(N),
$$
$R^N(d\varphi(e_j),d\varphi(e_k))d\varphi(e_k)$ is tangent to 
$\varphi_{\ast}TM$. By (2.16) of Lemma 2.2, 
\begin{equation}
\overline{\Delta}\tau(\varphi)=
\sum_{i,j=1}^mh(\tau(\varphi),(\widetilde{\nabla}_{e_i}d\varphi)(e_j))
(\widetilde{\nabla}_{e_i}d\varphi)(e_j).
\end{equation}
Furthermore, we have 
\begin{align}
{\mathcal R}(\tau(\varphi))
&=\sum_{i=1}^mR^N(\tau(\varphi),d\varphi(e_i))d\varphi(e_i)\nonumber\\
&=c\sum_{i=1}^m\{
h(d\varphi(e_i),d\varphi(e_i))\tau(\varphi)
-h(d\varphi(e_i),\tau(\varphi))d\varphi(e_i)
\}\nonumber\\
&=cm\tau(\varphi).
\end{align}
Then, $\varphi:\,(M,g)\rightarrow S^{m+1}(\frac{1}{\sqrt{c}})$ is biharmonic if and only if 
\begin{align}
\tau_2(\varphi)
&=\overline{\Delta}\tau(\varphi)-{\mathcal R}(\tau(\varphi))\nonumber\\
&=\sum_{i,j=1}^mh(\tau(\varphi),(\widetilde{\nabla}_{e_i}d\varphi)(e_j))
(\widetilde{\nabla}_{e_i}d\varphi)(e_j)-cm\tau(\varphi)\nonumber\\
&=0.
\end{align}
If we denote by $\xi$, the unit normal vector field to $\varphi(M)$, 
the second fundamental form $B(\varphi)$ is of the form
$B(\varphi)(e_i,e_j)=(\widetilde{\nabla}_{e_i}d\varphi)(e_j)=h_{ij}\xi$.  Then, we have 
$\tau(\varphi)=\sum_{i=1}^mB(\varphi)(e_i,e_i)=\sum_{i=1}^mh_{ii}\,\xi$ and $\Vert B(\varphi)\Vert^2=\sum_{i,j=1}^mh_{ij}h_{ij}$. 
Substituting these into (3.3), we have 
\begin{equation}
\tau_2(\varphi)=
\sum_{k=1}^mh_{kk}\left(
\sum_{i,j=1}^mh_{ij}h_{ij}-cm
\right)\xi=0, 
\end{equation}
That is, 
$\Vert B(\varphi)\Vert^2=cm$ 
since 
$\sum_{k=1}^mh_{kk}\not=0$. 
\end{pf}
\vskip0.6cm\par
Next, we prepare the necessary materials on isoparametric hypersurfaces $M$  in the unit sphere $S^n(1)$ following M\"unzner (\cite{M}) or Ozeki and Takeuchi (\cite{OT}). 
\par
Let $\varphi:\, (M,g)\rightarrow S^n(1)$ be an isometric immersion of 
$(M,g)$ into the unit sphere $S^n(1)$ and denote by $(N,h)$, the unit sphere $S^n(1)$ with the canonical metric. Assume that 
$\dim M=n-1$. 
The shape operator 
$A_{\xi}$ is a linear operator of $T_xM$ into itself defined by 
$$
g(A_{\xi}X,Y)=h(\varphi_{\ast}(\nabla_XY), \xi),\,X,Y\in {\frak X}(M),
$$ 
where $\xi$ is the unit normal vector field along $M$. 
The eigenvalues of $A_{\xi}$ are called the {\em principal curvatures}. $M$ is called {\em isoparametric} if all the principal curvatures 
are constant in $x\in M$.  
It is known that there exists a homogeneous polynomial $F$ on ${\Bbb R}^{n+1}$ of degree $g$ whose restriction to $S^n(1)$, 
denoted by $f$, called {\em isoparametric function}, 
$M$ is given by $M=f^{-1}(t)$ for some $t\in I=(-1,1)$. 
For each $t\in I$, 
$\xi_t=\frac{\nabla f}{\sqrt{g(\nabla f,\nabla f)}}$ is a smooth unit normal vector field along $M_t=f^{-1}(t)$, and all the distinct principal curvatures of $M_t$ with respect to $\xi_t$ 
are given as
$$
k_1(t)>k_2(t)>\cdots>k_{g(t)}(t)
$$
with their multiplicities $m_j(t)$ $(j=1,\cdots,g(t))$. And $g=g(t)$ is constant in $t$, and is should be 
$g=1,2,3,4, \,\,{\rm or}\,\, 6$.  Furthermore, it holds that 
\begin{align}
m_1(t)&=m_3(t)=\cdots=m_1,\nonumber\\
m_2(t)&=m_4(t)=\cdots=m_2,\nonumber\\
k_j(t)&=
\cot\left(
\frac{(j-1)\pi+\cos^{-1}t}{g}
\right)\quad (j=1,\cdots,g).
\end{align}
where $m_1$ and $m_2$ are constant in $t\in I$. 
We also have 
\begin{equation}
\Vert B(\varphi)\Vert ^2=\Vert A_x\Vert ^2=
\sum_{j=1}^{g(t)}m_j(t)k_j(t)^2.
\end{equation}
Indeed, if we denote by $\lambda_i$ $(i=1,\cdots,m$ $(m=\dim M)$, all the principal curvature counted with their multiplicities, 
we may choose orthonomal eigenvectors $\{X_i\}_{i=1}^m$ of $T_xM$ 
in such a way that $A_{\xi}X_i=\lambda_iX_i$ $(i=1,\cdots,m)$. 
Then, we have 
$h(B(X_i,X_j),\xi)=g(A_{\xi}(X_i),X_j)=\lambda_i\delta_{ij}$, 
and $\Vert B(X_i,X_j)\Vert^2=\lambda_i{}^2\delta_{ij}$. Thus, we have
\begin{prop} 
Let $\varphi:\,(M,g)\rightarrow S^n(1)$ be an isoparametric hypersurface 
in the unit sphere $S^n(1)$, $\dim M=n-1$. Then, 
\begin{equation}
\Vert B(\varphi)\Vert^2=\sum_{j=1}^m\lambda_j{}^2.
\end{equation}
\end{prop}
\begin{pf} Indeed, we have 
$$\Vert B(\varphi)\Vert^2=\Vert A_{\xi}\Vert ^2=\sum_{i,j=1}^m\Vert B(X_i,X_j)\Vert^2=\sum_{j=1}^m\lambda_j{}^2,$$
which is (3.7). 
\end{pf}
\vskip0.6cm\par
\section{Biharmonic isoparametric hypersurfaces}
Now, our main theorem in this section is 
\begin{th}
Let $\varphi:\, (M,g)\rightarrow S^n(1)$ be an isometric immersion $($$\dim M=n-1$$)$ which is isoparametric. Then, $(M,g)$ is biharmonic if and only if $(M,g)$ is one of the following:
\par
$({\rm i})$ \quad $M=S^{n-1}\left(\frac{1}{\sqrt{2}}
\right) \subset S^n(1)$, \quad $($a small sphere$)$
\par
$({\rm ii})$ \quad 
$M=S^{n-p}\left(\frac{1}{\sqrt{2}}
\right)\times S^{p-1}\left(\frac{1}{\sqrt{2}}
\right)\subset S^n(1)$, with $n-p\not=p-1$
\par\qquad \quad
$($the Clifford torus$)$, or
\par
$({\rm iii})$ \quad $\varphi:\,(M,g)\rightarrow S^n(1)$ is harmonic, i.e., minimal. 
\end{th}
\begin{pf}
The proof is divided into the cases $g=1,2,3,4, {\rm or} \,6$. 
It is known that for the cases $g=1,2$, all the $(M,g)$ are homogeneous, and are classified into two cases. 
For $g=3,4$ or $6$, we will show there are no nonharmonic biharmonic isoparametric hypersurfaces in the unit sphere. 
\par
\underline{Case 1: $g=1$}. \quad In this case, $m_1=m_2=n-1$ and 
$k_1(t)=\cot x$, $x=\cos^{-1}t$ with $0<x<\pi$, $-1<t<1$. 
Then, we have immediately: 
\par\qquad\qquad
minimal $\Longleftrightarrow\cot x=0\Longleftrightarrow t=0$ 
(a great sphere). 
\par\noindent
Furthermore, we have: 
\par\qquad
biharmonic and nonminimal
$\Longleftrightarrow (n-1)\cot^2x=n-1$
\par\qquad\qquad\quad\qquad\qquad\qquad\qquad\quad\,\,
$\Longleftrightarrow t=\pm\frac{1}{\sqrt{2}}$\,(a small sphere).
\par
\underline{Case 2: $g=2$}.  \quad 
In this case, $m_1=p-1$, $m_2=n-p$ with 
$(2\leq p\leq\left[\frac{n+1}{2}\right])$. 
Then, we have immediately, 
\par
\begin{align}
{\rm minimal}&\Longleftrightarrow (p-1)\cot\left(\frac{x}{2}\right)+(n-p)\cot\left(
\frac{\pi+x}{2}
\right)=0\nonumber\\
&\Longleftrightarrow\cos^2\left(\frac{x}{2}
\right)=\frac{n-p}{n-1}\nonumber\\
&\Longleftrightarrow t=\frac{n+1-2p}{n-1}\nonumber
\end{align}
with $x=\cos^{-1}t$ $(t\in (-1,1))$.  
On the other hand, by Proposition 3.1, 
\begin{align}
{\rm biharmonic}
&\Longleftrightarrow 
(p-1)\cot^2\left(\frac{x}{2}
\right)+(n-p)\cot^2\left(\frac{\pi+x}{2}=n-1\right)\nonumber\\
&\Longleftrightarrow 
t=0,\,\frac{n+1-2p}{n-1}\nonumber
\end{align}
with $x=\cos^{-1}t$ $(t\in (-1,1))$.  Thus, 
\par
biharmonic  and nonminimal
$\Longleftrightarrow 
t=0$, 
\par\qquad 
$k_1(0)=1$ $(m_1=p-1)$,\,\,$k_2(0)=-1$ $(m_2=n-p)$
\par\qquad
$p-1\not= n-p$. 
\par
\underline{Case 3: $g=3$}. 
In this case, all the isoparametric hypersurfaces are classified into four cases, and $m_1=m_2$ are $1,2,4$ or, $8$, and $\dim M$ is 
$3, 6, 12$ or $24$, respectively. By Proposition 3.1, it suffices to show in the case $\dim M=3$, 
\begin{equation}
\cot^2\left(\frac{x}{3}\right)
+\cot^2\left(\frac{\pi+x}{3}\right)
+\cot^2\left(\frac{2\pi+x}{3}\right)
\geq 6>3\quad (0<x<\pi).
\end{equation} 
To prove (4.3), 
we only see the LHS of (4.3) coincides with 
\begin{equation}
\cot^2\left(\frac{x}{3}\right)
+\left(
\frac{\cot\frac{x}{3}-\sqrt{3}}{\sqrt{3}\cot\frac{x}{3}+1}
\right)^2
+\left(
\frac{\cot\frac{x}{3}+\sqrt{3}}{-\sqrt{3}\cot\frac{x}{3}+1}
\right)^2,
\end{equation}
which is bigger than or equal to $6$ when $0<x<\pi$. 
Remark that $0<\cot\frac{x}{3}<\frac{1}{\sqrt{3}}$ $(0<x<\pi)$. 
And the arguments go the same way as $\dim M=6, 12, 24$. 
Thus, due to Proposition 3.1 and Theorem 3.1, 
 there are no nonminimal biharmonic hypersurfaces in this case. 
\par
\underline{Case 4: $g=4$}. 
In this case, we have 
\begin{align}
\Vert B(\varphi)\Vert^2
&=m_1(t)\cot^2\left(\frac{x}{4}\right)
+m_2(t)\cot^2\left(\frac{\pi+x}{4}\right)\nonumber\\
&\quad
+m_1(t)\cot^2\left(\frac{2\pi+x}{4}\right)
+m_2(t)\cot^2\left(\frac{3\pi+x}{4}\right)\nonumber\\
&=m_1(t)\left\{
\cot^2\left(\frac{x}{4}\right)+\frac{1}{\cot^2\left(\frac{x}{4}\right)}
\right\}\nonumber\\
&+m_2(t)
\left\{
\left(
\frac{\cot\left(\frac{x}{4}\right)-1}{\cot\left(\frac{x}{4}\right)+1}
\right)^2
+
\left(\frac{\cot\left(\frac{x}{4}\right)+1}{\cot\left(\frac{x}{4}\right)-1}\right)^2
\right\}\nonumber\\
&\geq 2m_1(t)+2m_2(t)=\dim M, 
\end{align}
and equality holds if and only if 
\begin{equation}
\left\{
\begin{aligned}
&\cot^2\left(\frac{x}{4}\right)=\frac{1}{\cot^{2}\left(\frac{x}{4}\right)},
\\
&\left(
\frac{\cot\left(\frac{x}{4}\right)-1}{\cot\left(\frac{x}{4}\right)+1}
\right)^2
=
\left(\frac{\cot\left(\frac{x}{4}\right)+1}{\cot\left(\frac{x}{4}\right)-1}\right)^2,
\end{aligned}
\right.
\end{equation}
because, for all $a>0$ and $b>0$, 
$\frac{a+b}{2}\geq\sqrt{ab}$
and equality holds if and only if $a=b$. 
But, it is impossible that (4.6) holds. Thus, we have 
$\Vert B(\varphi)\Vert^2>\dim M$. In this case, due to Proposition 3.1 and Theorem 3.1, there are no nonharmonic biharmonic immersions $\varphi$. 
\par
\underline{Case 5: $g=6$}. In this case, we have 
\begin{align}
&\Vert B(\varphi)\Vert^2
=
m_1(t)\cot^2\left(\frac{x}{6}\right)
+m_2(t)\cot^2\left(\frac{\pi+x}{6}\right)\nonumber\\
&\qquad\qquad
+m_1(t)\cot^2\left(\frac{2\pi+x}{6}\right)
+m_2(t)\cot^2\left(\frac{3\pi+x}{6}\right)\nonumber\\
&\qquad\qquad
+m_1(t)\cot^2\left(\frac{4\pi+x}{6}\right)
+m_2(t)\cot^2\left(\frac{5\pi+x}{6}\right)\nonumber\\
&=m_1(t)\left\{
\cot^2\left(\frac{x}{6}\right)
+\left(
\frac{\cot\frac{x}{6}-\sqrt{3}}{\sqrt{3}\cot\frac{x}{6}+1}
\right)^2
+\left(
\frac{\cot\frac{x}{6}+\sqrt{3}}{-\sqrt{3}\cot\frac{x}{6}+1}
\right)^2
\right\}\nonumber\\
&\,\,+
m_2(t)
\left\{
\frac{1}{\cot^2\left(\frac{x}{6}\right)}
+
\left(
\frac{\sqrt{3}\cot\frac{x}{6}-1}{\cot\frac{x}{6}+\sqrt{3}}
\right)^2
+\left(
\frac{\sqrt{3}\cot\frac{x}{6}+1}{-\cot\frac{x}{6}+\sqrt{3}}
\right)^2
\right\}.
\end{align}
Here, we denote by $f(y)$, the bracket of the first term of the RHS of (4.7), 
where $y=\cot\frac{x}{6}>\sqrt{3}$ $(0<x<\pi)$. 
Then, we have 
$\frac{df}{dy}>0$ and $\lim_{y\rightarrow \sqrt{3}}f(y)=6$. 
And we denote by $g(y)$, the bracket of the second term of the RHS of (4.7), 
where $y=\cot\frac{x}{6}>\sqrt{3}$ $(0<x<\pi)$. 
Then, we have 
$\frac{dg}{dy}<0$ and $\lim_{y\rightarrow \infty}g(y)=6$. 
Therefore, we have
\begin{equation}
\Vert B(\varphi)\Vert^2
\geq 6(m_1(t)+m_2(t))>3(m_1(t)+m_2(t))=\dim M.
\end{equation}
Thus, due to Proposition 3.1 and Theorem 3.1, there are also no nonharmonic biharmonic immersions $\varphi$ 
in this case. 
\end{pf} 
\vskip0.6cm\par
\section{Biharmonic maps into the complex projective space}
In the following two sections, we show classification of all  homogeneous real hypersurfaces in the complex $n$-dimensional projective space ${\Bbb C}P^n(c)$ 
with positive constant holomorphic sectional curvature $c>0$
which are {\em biharmonic}. 
To do it, we need first the following theorem analogue to Theorem 3.1 which charcterizes the biharmonic maps. 
\begin{th} Let $(M,g)$ be a real $(2n-1)$-dimensional compact Riemannian manifold, and
$\varphi:\,(M,g)\rightarrow {\Bbb C}P^n(c)$ be an isometric immersion with non-zero constant mean curvature. Then, the necessary and sufficient condition for $\varphi$ to be biharmonic is 
\begin{equation}
\Vert B(\varphi)\Vert ^2=\frac{n+1}{2}c. 
\end{equation}
\end{th}
\begin{pf}
By Lemma 2.3, the mean curvature vector of $\varphi$ is parallel with respect to $\nabla^{\perp}$, so we may apply Lemmas 2.1 and 2.2 in this case. 
Let us recall the fact that the curvature tensor of $(N,h)={\Bbb C}P^n(c)$ is given by 
\begin{align}
R^N(U,V)W&=\frac{c}{4}
\big\{
h(V,W)U-h(U,W)V\nonumber\\
&\qquad
+h(JV,W)JU-h(JU,W)JV+2h(U,JV)JW\big\},\nonumber
\end{align}
where $J$ is the adapted almost complex tensor, and 
$U$,$V$ and $W$ are vector fields on 
${\Bbb C}P^n(c)$. 
Then, we have 
\begin{align}
R^N(d\varphi(e_j),&d\varphi(e_k))d\varphi(e_k)=
\frac{c}{4}\big\{d\varphi(e_j)-\delta_{jk}\,d\varphi(e_k)\nonumber\\
&\qquad\quad+3h(d\varphi(e_j),Jd\varphi(e_k))\,Jd\varphi(e_k)\big\}. 
\end{align}
Then, we have 
\begin{equation}
\sum_{j,k=1}^mh\big(\tau(\varphi),
R^N(d\varphi(e_j),d\varphi(e_k))d\varphi(e_k)\big)\,
d\varphi(e_j)=0.
\end{equation}
Because the LHS of (5.3) coincides with 
\begin{align}
\frac{3c}{4}&\sum_{j,k=1}^mh\big(d\varphi(e_j),Jd\varphi(e_k)\big)\,
h\big(\tau(\varphi),Jd\varphi(e_k)\big)\,d\varphi(e_j)\nonumber\\
&\quad=
\frac{3c}{4}\sum_{j,k=1}^mh\big(Jd\varphi(e_j),d\varphi(e_k)\big)\,
h\big(J\tau(\varphi),d\varphi(e_k)\big)\,d\varphi(e_j)\nonumber\\
&\quad=
\frac{3c}{4}\sum_{j=1}^mh\big(Jd\varphi(e_j),
\sum_{k=1}^m
h\big(J\tau(\varphi),d\varphi(e_k)\big)d\varphi(e_k)\big)\,
d\varphi(e_j)\nonumber\\
&\quad=
\frac{3c}{4}\sum_{j=1}^mh\big(Jd\varphi(e_j),J\tau(\varphi)\big)\,
d\varphi(e_j)\nonumber\\
&\quad=
\frac{3c}{4}\sum_{j=1}^mh(d\varphi(e_j),\tau(\varphi))\,
d\varphi(e_j)=0.
\end{align}
Here the third equality follows from that
$J\tau(\varphi)\in \Gamma(\varphi_{\ast}TM)$ 
which is due to 
$h(J\tau(\varphi),\tau(\varphi))=0$, $0\not=\tau(\varphi)\in T^{\perp}M$  and $\dim M=2n-1$.  
Since $\{d\varphi(e_k)\}_{k=1}^m$ is an orthonormal basis 
of $\varphi_{\ast}(T_xM)$ at each $x\in M$, 
$J\tau(\varphi)=\sum_{k=1}^m
h(J\tau(\varphi),d\varphi(e_k))d\varphi(e_k)$. 
\par
By (2.16) in Lemma 2.2, we have 
\begin{equation}
\overline{\Delta}\tau(\varphi)
=\sum_{i,j=1}^m h\big(\tau(\varphi),(\widetilde{\nabla}_{e_i}d\varphi)(e_j)\big)\,
(\widetilde{\nabla}_{e_i}d\varphi)(e_j).
\end{equation}
Furthermore, we have 
\begin{equation}
{\mathcal R}(\tau(\varphi))
=\frac{c}{4}(m+3)\tau(\varphi).
\end{equation}
Because 
the LHS of (5.6) is equal to 
\begin{align}
\sum_{k=1}^mR^N(\tau(\varphi),d\varphi(e_k))d\varphi(e_k)
&=\frac{c}{4}\big\{m\tau(\varphi)\nonumber\\
&\qquad
-3\sum_{k=1}^m
h(J\tau(\varphi),d\varphi(e_k))\,Jd\varphi(e_k)\big\}\nonumber\\
&=\frac{c}{4}\{m\tau(\varphi)-3J(J\tau(\varphi))\}\nonumber\\
&=\frac{c}{4}(m+3)\tau(\varphi).
\end{align}
\par
Now the sufficient and necessary condition for $\varphi$ to be biharmonic 
is that 
\begin{equation}
\tau_2(\varphi)=\overline{\Delta}\tau(\varphi)
-{\mathcal R}
(\tau(\varphi))=0
\end{equation}
which is equivalent to 
\begin{equation}
\sum_{i,j=1}^m
h(\tau(\varphi),(\widetilde{\nabla}_{e_i}d\varphi)(e_j))(\widetilde{\nabla}_{e_i}d\varphi)(e_j)
-\frac{c}{4}(m+3)\tau(\varphi)=0.
\end{equation}
Here, we may denote as
\begin{align}
&B(\varphi)(e_i,e_j)=(\widetilde{\nabla}_{e_i}d\varphi)(e_j)=h_{ij}\,\xi\nonumber\\
&\tau(\varphi)=\sum_{k=1}^m(\widetilde{\nabla}_{e_k}d\varphi)(e_k)
=\sum_{k=1}^mh_{kk}\xi,
\end{align}
where $\xi$ is the unit normal vector field along $\varphi(M)$. 
Thus, the LHS of (5.9) coincides with 
\begin{align}
\sum_{i,j,k=1}^m&h_{kk}h_{ij}h_{ij}-\frac{c}{4}(m+3)\sum_{k=1}^mh_{kk}\nonumber\\
&=
\bigg(\sum_{k=1}^mh_{kk}\bigg)
\left\{\sum_{i,j=1}^mh_{ij}h_{ij}-\frac{c}{4}(m+3)
\right\}\nonumber\\
&=\Vert\tau(\varphi)\Vert^2
\,\left\{
\Vert B(\varphi)\Vert^2-\frac{c}{2}(n+1)
\right\},
\end{align} 
which yields the desired (5.1) due to the assumption that 
$\Vert \tau(\varphi)\Vert$ is a non-zero constant. 
\end{pf}
\vskip0.6cm\par
\section{Biharmonic Homogeneous real hypersurfaces in the complex projective space}
In this section, we classify all the {\em biharmonic} homogeneous real hypersurfaces in the complex projective space 
${\Bbb C}P^n(c)$. 
\par
First, let us recall the classification theorem of all the homogeneous real hypersurfaces in ${\Bbb C}P^n(c)$ due to R. Takagi 
(cf. \cite{T1}) based on a work by W.Y. Hsiang and H.B. Lawson (\cite{HL}). Let $U/K$ be a symmetric space of rank two of compact type, and ${\frak u}={\frak k}\oplus {\frak p}$, the Cartan decomposition of the Lie algebra $\frak u$ of $U$, and the Lie subalgebra 
$\frak k$ corresponding to $K$. 
Let $\langle X,\,Y\rangle=-B(X,Y)$ 
$(X,Y\in {\frak p})$ be the inner product on ${\frak p}$, 
$\Vert X\Vert^2=\langle X,X\rangle$, and 
$S:=\{X\in {\frak p};\,\Vert X\Vert =1\}$, the unit sphere in the Euclidean space $({\frak p}, \langle\,,\,\rangle)$, where 
$B$ is the Killing form of ${\frak u}$. Consider the adjoint action of $K$ on $\frak p$. Then, the orbit $\hat{M}={\rm Ad}(K)A$ 
through any regular element $A\in {\frak p}$ with $\Vert A\Vert=1$ gives a homogeneous hypersurface in the unit sphere $S$.  Conversely, any homogeneous hypersurface in $S$ can be obtained in this way (\cite{HL}). 
\par
Let us take as $U/K$, a {\em Hermitian} symmetric space of compact type of rank two of complex dimension $(n+1)$, and identify 
${\frak p}$ with ${\Bbb C}^{n+1}$.  Then, the adjoint orbit $\hat{M}={\rm Ad}(K)A$ of $K$ through any regular element $A$ in ${\frak p}$ 
is again a homogeneous hypersurface in the unit sphere $S$. 
Let $\pi:\,{\Bbb C}^{n+1}-\{{\bf 0}\}={\frak p}-\{{\bf 0}\}\rightarrow {\Bbb C}P^n$ be the natural projection. Then, the projection induces the Hopf fibration of $S$ onto ${\Bbb C}P^n$, denoted also by $\pi$, and 
$\varphi:\,M:=\pi(\hat{M})\hookrightarrow {\Bbb C}P^n$ gives a homogeneous real hypersurface in the complex projective space ${\Bbb C}P^n(4)$ with constant holomorphic sectional curvature $4$.  Conversely, any homogeneous real hypersurface $M$ in ${\Bbb C}P^n(4)$ is given in this way (\cite{T1}). 
Furthermore, all such hypersurfaces are classified into the following five types: 
\vskip0.3cm\par
(1) $A$-type: \par
${\frak u}={\frak su}(p+2)\oplus{\frak su}(q+2)$, 
${\frak k}={\frak s}({\frak u}(p+1)+{\frak u}(1))
\oplus{\frak s}({\frak u}(q+1)+{\frak u}(1))$, where 
$0\leq p\leq q$, $0<q$,  $p+q=n-1$, and $\dim M=2n-1$. 
\par
(2) $B$-type: \,
${\frak u}={\frak o}(m+2)$, ${\frak k}={\frak o}(m)\oplus {\Bbb R}$, where $3\leq m$, $\dim M=2m-3$. 
\par
(3) $C$-type:\quad ${\frak u}={\frak su}(m+2)$, ${\frak k}={\frak s}({\frak o}(m)+ {\frak o}(2))$ , where $3\leq m$, and $\dim M=4m-3$. 
\par
(4) $D$-type:\quad ${\frak o}(10)$, ${\frak u}(5)$, and 
$\dim M=17$. 
\par
(5) $E$-type:\quad 
${\frak u}={\frak e}_6$, ${\frak k}={\frak o}(10)\oplus{\Bbb R}$, and 
$\dim M=29$. 
\vskip0.3cm\par
He also gave (\cite{T2}, \cite{T3}) lists of the principal curvatures and their multiplicities of these $M$ as follows: 
\vskip0.3cm\par
(1) $A$-type:  Assume that 
$$
U/K=\frac{SU(p+2)\times SU(q+2)}{S(U(p+1)\times U(1))\times S(U(q+1)\times U(1)},
$$
then, the adjoint orbit of $K$, ${\rm Ad}(K)A$ is given by the Riemannian product of two odd dimensional spheres,
\begin{equation}\hat{M}=\hat{M}_{p,q}=S^{2p+1}(\cos u)\times S^{2q+1}(\sin u)\subset S^{2n+1},
\end{equation}
where $0<u<\frac{\pi}{2}$. 
The projection $M_{p,q}(u):=\pi(\hat{M}_{p,q}(u))$ is 
a homogeneous real hypersurface of ${\Bbb C}P^n(4)$. 
The principal curvatures of $M_{p,q}$ with $0\leq p\leq q$, $0<q$, 
are given as 
\begin{equation}
\left\{
\begin{aligned}
\lambda_1&=-\tan u \quad ({\rm with\,\, multiplicity}\,\, m_1=2p),
\\
&\qquad\qquad(m_1=0\,\,{\rm if}\,\,p=0),\\
\lambda_2&=\cot u\quad ({\rm with\,\,multiplicity}\,\,m_2=2q),\\
\lambda_3&=2\cot(2u)\quad ({\rm with\,\,multiplicity}\,\,m_3=1).
\end{aligned}
\right.
\end{equation}
\par
Thus, the mean curvature $H$ of $M_{p,q}(u)$ is given by 
\begin{align}
H&=\frac{1}{2n-1}\{
2q\cot u-2p\tan u+2\cot (2u)\}
\nonumber\\
&=\frac{1}{2n-1}\{
(2q+1)\cot u-(2p+1)\tan u\}.
\end{align}
The constant 
$\Vert B(\varphi)\Vert^2$
which is the sum of all the square of principal curvatures 
with their multiplicities, is given by 
\begin{align}
\Vert B(\varphi)\Vert ^2
&=2q\cot^2u+2p\tan^2u+4\cot^2(2u)\nonumber\\
&=(2q+1)\cot^2u+(2p+1)\tan^2u-2.
\end{align}
\vskip0.3cm\par
(2) $B$-type: \quad Assume that 
$U/K=SO(m+2)/(SO(m)\times SO(2))$, 
$(m:=n+1)$, and then, 
the adjoint orbit of $K$, ${\rm Ad}(K)A$ is given by 
$$\hat{M}=\{ SO(n+1)\times SO(2)\}/\{SO(n-1)\times {\Bbb Z}_2\}\subset S^{2n+1}.$$
The real hypersurface $\varphi:\,M\hookrightarrow {\Bbb C}P^n$ is a tube over a complex quadric with radius 
$\frac{\pi}{4}-u$ $(0<u<\frac{\pi}{4})$ 
or a tube over a totally geodesic 
real projective space ${\Bbb R}P^n$ with radius $u$ 
$(0<u<\frac{\pi}{4})$. The principal curvatures of $M$ are given as 
\begin{equation}
\left\{
\begin{aligned}
\lambda_1&=-\cot u \quad ({\rm with\,\, multiplicity}\,\, m_1=n-1),
\\
\lambda_2&=\tan u\quad ({\rm with\,\,multiplicity}\,\,m_2=n-1),\\
\lambda_3&=2\tan(2u)\quad ({\rm with\,\,multiplicity}\,\,m_3=1).
\end{aligned}
\right.
\end{equation}
\par
Thus, the mean curvature of $M$ is given by 
\begin{align}
H&=\frac{1}{2n-1}\{-(n-1)\cot u+(n-1)\tan u+2\cot(2u)\}\nonumber\\
&=-\frac{1}{2n-1}\,\cdot\,\frac{(n-1)t^4-2(n+1)t^2+n-1}{t(t^2-1)},
\end{align}
where $t=\cot u$. 
The constant $\Vert B(\varphi)\Vert ^2$ is given by 
\begin{align}
\Vert B(\varphi)\Vert^2&=(n-1)\cot^2u+{n-1}\tan^2u+4\tan^2(2u)\nonumber\\
&=(n-1)t^2+\frac{n-1}{t^2}+\frac{16t^2}{(t^2-1)^2}\nonumber\\
&=\frac{(n-1)(X-1)^2(X^2+1)+16X^2}{X(X-1)^2},
\end{align}
where 
$X:=t^2$. 
\vskip0.3cm\par
(3) $C$-type: \quad Assume that 
$U/K=SU(m+2)/S(U(m)\times U(2))$, $(n=2m+1)$, and then, 
the adjoint orbit of $K$, ${\rm Ad}(K)A$ is given by 
$$
\hat{M}=S(U(m)\times U(2))/(T^2\times SU(m-2))\subset S^{2n+1}.
$$
The real hypersurface $\varphi:\,M\hookrightarrow {\Bbb C}P^n$ is 
a tube over the Segre imbeding of ${\Bbb C}^1\times {\Bbb C}P^m$ with radius $u$ $(0<u<\frac{\pi}{4})$. The principal curvatures of $M$ are given by 
\begin{equation}
\left\{
\begin{aligned}
\lambda_1&=-\cot u \quad ({\rm with\,\, multiplicity}\,\, m_1=n-3),
\\
\lambda_2&=\cot \left(\frac{\pi}{4}-u\right)
\quad ({\rm with\,\,multiplicity}\,\,m_2=2),\\
\lambda_3&=\cot\left(\frac{\pi}{2}-u\right) 
\quad ({\rm with\,\,multiplicity}\,\,m_3=n-3),\\
\lambda_4&=\cot\left(\frac{3\pi}{4}-u\right)
\quad ({\rm with\,\,multiplicity}\,\,m_4=2),\\
\lambda_5&=-2\tan(2u)\quad ({\rm with\,\,multiplicity}\,\,m_5=1).
\end{aligned}
\right.
\end{equation}
\par
Then, 
$$
\lambda_1=-t,\,\lambda_2=\frac{t+1}{t-1},\,\lambda_3=\frac{1}{t},\,
\lambda_4=-\frac{t-1}{t+1},\,\lambda_5=-t+\frac{1}{t},
$$
where $t=\cot u$. 
The mean curvature of $M$ is given by 
\begin{align}
H&=
\frac{1}{2n-1}\left\{
(n-3)(-t)+2\frac{t+1}{t-1}+(n-3)\,\frac{1}{t}-2\frac{t-1}{t+1}-t+\frac{1}{t}
\right\}
\nonumber\\
&=-\frac{(n-2)t^4-2(n+2)t^2+n-2}{t(t^2-1)}.
\end{align}
The constant $\Vert B(\varphi)\Vert ^2$ is given by 
\begin{align}
\Vert B(\varphi)\Vert^2
&=(n-3)t^2+2\left(\frac{t+1}{t-1}\right)^2
+(n-3)\frac{1}{t^2}
\nonumber\\
&\quad+2\left(\frac{t-1}{t+1}\right)^2
+\left(-t+\frac{1}{t}\right)^2
\nonumber\\
&=\frac{C(X)}{X(X-1)^2},
\end{align}
where 
\begin{align}
C(X)&:=(n-2)X^2(X-1)^2+(n-2)(X-1)^2\nonumber\\
&\quad+4X(X^2+6X+1)-2X(X-1)^2,
\end{align}
and 
$X:=t^2$. 
\vskip0.3cm\par
(4) $D$-type: \quad Assume that 
$U/K=O(10)/U(5)$, and then, 
the adjoint orbit of $K$, ${\rm Ad}(K)A$ is given by 
$$
\hat{M}=U(5)/(SU(2)\times SU(2)\times U(1))\subset S^{19}.
$$
The real hypersurface $\varphi:\,M\hookrightarrow {\Bbb C}P^9$ is a tube over the Pl\"ucker 
 imbeding of ${\rm Gr}_2({\Bbb C}^5)$ with radius 
 $u$ $(0<u<\frac{\pi}{4})$. 
 The principal curvatures of $M$ are given by 
 \begin{equation}
\left\{
\begin{aligned}
\lambda_1&=-\cot u \quad ({\rm with\,\, multiplicity}\,\, m_1=4),
\\
\lambda_2&=\cot \left(\frac{\pi}{4}-u\right)
\quad ({\rm with\,\,multiplicity}\,\,m_2=4),\\
\lambda_3&=\cot\left(\frac{\pi}{2}-u\right) 
\quad ({\rm with\,\,multiplicity}\,\,m_3=4),\\
\lambda_4&=\cot\left(\frac{3\pi}{4}-u\right)
\quad ({\rm with\,\,multiplicity}\,\,m_4=4),\\
\lambda_5&=-2\tan(2u)\quad ({\rm with\,\,multiplicity}\,\,m_5=1).
\end{aligned}
\right.
\end{equation}
\par
Then, 
$$
\lambda_1=-t,\,\lambda_2=\frac{t+1}{t-1},\,\lambda_3=\frac{1}{t},\,
\lambda_4=-\frac{t-1}{t+1},\,\lambda_5=-t+\frac{1}{t},
$$
where $t=\cot u$. 
The mean curvature of $M$ is given by 
\begin{align}
H&=
\frac{1}{17}\left\{
4(-t)+4\frac{t+1}{t-1}+4\,\frac{1}{t}-4\frac{t-1}{t+1}-t+\frac{1}{t}
\right\}
\nonumber\\
&=-\frac{5t^4-26t^2+5}{17t(t^2-1)}=-\frac{(5t^2-1)(t^2-5)}{17t(t^2-1)}.
\end{align}
The constant $\Vert B(\varphi)\Vert ^2$ is given by 
\begin{align}
\Vert B(\varphi)\Vert^2
&=4t^2+4\left(\frac{t+1}{t-1}\right)^2
+4\frac{1}{t^2}
\nonumber\\
&\quad+4\left(\frac{t-1}{t+1}\right)^2
+\left(-t+\frac{1}{t}\right)^2
\nonumber\\
&=\frac{D(X)}{X(X-1)^2},
\end{align}
where 
\begin{equation}
D(X):=11X^3+63X^2+X+5,
\end{equation}
and 
$X:=t^2$. 
\vskip0.3cm\par
(5) $E$-type: \quad Assume that 
$U/K=E_6/({\rm Spin}(10)\times U(1)$, and then, 
the adjoint orbit of $K$, ${\rm Ad}(K)A$ is given by 
$$
\hat{M}=({\rm Spin}(10)\times U(1))/(SU(4)\times U(1))\subset S^{31}.
$$
The real hypersurface $\varphi:\,M\hookrightarrow {\Bbb C}P^{15}$ is a tube over the canonical
 imbeding of $SO(10)/U(5)\subset {\Bbb C}P^{15}$ with radius 
 $u$ $(0<u<\frac{\pi}{4})$. 
 The principal curvatures of $M$ are given by 
 \begin{equation}
\left\{
\begin{aligned}
\lambda_1&=-\cot u \quad ({\rm with\,\, multiplicity}\,\, m_1=8),
\\
\lambda_2&=\cot \left(\frac{\pi}{4}-u\right)
\quad ({\rm with\,\,multiplicity}\,\,m_2=6),\\
\lambda_3&=\cot\left(\frac{\pi}{2}-u\right) 
\quad ({\rm with\,\,multiplicity}\,\,m_3=8),\\
\lambda_4&=\cot\left(\frac{3\pi}{4}-u\right)
\quad ({\rm with\,\,multiplicity}\,\,m_4=6),\\
\lambda_5&=-2\tan(2u)\quad ({\rm with\,\,multiplicity}\,\,m_5=1).
\end{aligned}
\right.
\end{equation}
\par
Then, 
$$
\lambda_1=-t,\,\lambda_2=\frac{t+1}{t-1},\,\lambda_3=\frac{1}{t},\,
\lambda_4=-\frac{t-1}{t+1},\,\lambda_5=-t+\frac{1}{t},
$$
where $t=\cot u$. 
The mean curvature of $M$ is given by 
\begin{align}
H&=
\frac{1}{29}\left\{
8(-t)+6\frac{t+1}{t-1}+8\,\frac{1}{t}-6\frac{t-1}{t+1}-t+\frac{1}{t}
\right\}
\nonumber\\
&=-\frac{9t^4-42t^2+9}{29t(t^2-1)}.
\end{align}
The constant $\Vert B(\varphi)\Vert ^2$ is given by 
\begin{align}
\Vert B(\varphi)\Vert^2
&=8t^2+6\left(\frac{t+1}{t-1}\right)^2
+8\frac{1}{t^2}
\nonumber\\
&\quad+6\left(\frac{t-1}{t+1}\right)^2
+\left(-t+\frac{1}{t}\right)^2
\nonumber\\
&=\frac{E(X)}{X(X-1)^2}-2,
\end{align}
where 
\begin{equation}
E(X):=21X^3+99X^2-9X+9,
\end{equation}
and 
$X:=t^2$.
\vskip0.6cm\par
Now we want to show the following:
\begin{th}
Let $M$ be any homogeneous real hypersurface in ${\Bbb C}P^n(4)$, so that $M$ is a tube of $A\sim E$  type. 
\par
$({\rm I})$ Then, for each type, there is a unique $u$ with $0<u< \frac{\pi}{4}$ in such a way that 
$M$ is a tube of radius $u$ and is minimal. 
\par
$({\rm II})$ Assume that $M$ is a biharmonic but not minimal. 
Then, $M$ is one of type $A$, $D$ or $E$. More precisely,  
\par
$(1)$ in the case of $A$-type, 
$M$ is a tube $M_{p,q}(u)$ 
of ${\Bbb C}P^p\subset {\Bbb C}P^n$  
$(p\geq 0$ and $q=(n-1)-p)$ of radius $u$ with $0<u<\frac{\pi}{2}$ of which 
$t=\cot u$ is a solution of the equation 
\begin{equation}\cot u
=\left\{
\frac{p+q+3\pm\sqrt{(p-q)^2+4(p+q+2)}}{1+2q}
\right\}^{1/2}.
\end{equation}
\par
$(2)$ In the case of $D$-type,
$M$ is a tube of the Pl\"ucker imbedding ${\rm Gr}_2({\Bbb C}^5)\subset {\Bbb C}P^9$ of radius $u$ with $0<u<\frac{\pi}{4}$ of which 
$t=\cot u$ is a unique solution of the equation
\begin{equation}
41t^6+43t^4+41t^2-15=0.
\end{equation}
I.e., $u=1.0917\cdots$. 
\par
$(3)$ In the case of $E$-type, 
$M$ is a tube of 
the imbedding $SO(10)/U(5)\subset {\Bbb C}P^{15}$ of radius $u$ with $0<u<\frac{\pi}{4}$ of which $t=\cot u$ is a unique solution of the equation
\begin{equation}
13t^6-107t^4+43t^2-9=0.
\end{equation}
I.e., $u=0.343448\cdots$. 
\end{th}
\begin{pf}
We give a proof case by case. 
\par
\underline{Case (1) $A$-Type}: By (6.3), $\varphi:\,M_{p,q}(u)\hookrightarrow {\Bbb C}P^n(4)$ is harmonic if and only if 
\begin{equation}
t:=\cot u=\left\{
\frac{2p+1}{2q+1}
\right\}^{1/2}.
\end{equation}
On the other hand, by Theorem 5.1 and (6.4), 
$\varphi:\,M_{p,q}\hookrightarrow {\Bbb C}P^n(4)$ 
is non-harmonic and biharmonic if and only if 
$t=\cot u$ must satisfy 
\begin{equation} 
(2q+1)\cot^4u-2(p+q+3)\cot^2u+2p+1=0,
\end{equation}
so that 
\begin{equation}
t=\cot u=\left\{
\frac{p+q+3\pm\sqrt{(p-q)^2+4(p+q+2)}}{2q+1}
\right\}^{1/2}
\end{equation}
since 
$p+q+3\pm\sqrt{(p-q)^2+4(p+q+2)}$ is positive but does never be $2p+1$. 
\par
 \underline{Case (2) $B$-Type}: 
 By (6.6), $\varphi:\,M\hookrightarrow {\Bbb C}P^n(4)$ is harmonic if and only if 
$t=\cot u$  $(0<u<\frac{\pi}{4})$ 
must satisfy 
\begin{equation}
(n-1)t^4-2(n+1)t^2+n-1=0,
\end{equation}
which is equivalent to that 
\begin{equation}
t=\cot u=\left\{
\frac{n+1\pm2\sqrt{n}}{n-1}
\right\}^{1/2}
=\frac{\sqrt{n}\pm1}{\sqrt{n-1}}.
\end{equation}
On the other hand, by Theorem 5.1 and (6.7), 
$\varphi:\,M\hookrightarrow {\Bbb C}P^n(4)$ is non-harmonic but biharmonic if and only if 
\begin{align}
f(X):&=(n-1)(X-1)^2(X^2+1)+16X^2-2(n+1)X(X-1)^2\nonumber\\
&=0,
\end{align}
where $X:=t^2$. But, 
$f(X)>0$ for all $0<X<\infty$. Indeed, 
$(1)$ we have 
$$
f(X)=(n-1)(X-1)^2\left\{
X^2-2\frac{n+1}{n-1}X+1
\right\}+16X^2,
$$
which is positive when either
$X\geq 4$ and $n\geq 3$ or $X\leq 0.2679$ and $n\geq 3$. 
Furthermore, (2) 
we have 
$$
f(X)=(n-1)(X-1)^4+4X\left(4X-(X-1)^2\right),
$$
and 
$4X-(X-1)^2>0$ if 
$0.171573=3-2\sqrt{2}<X<3+2\sqrt{2}=5.82843$. 
So we have, 
$f(X)>0$ when $0.172<X<5.82$, 
Thus, by (1) and (2), $f(X)>0\,(0<X<\infty)$ when $n\geq 3$. 
In the case $n=2$, $f(X)=X^4-8X^3+30X^2-8X+1>0$ 
on $(0,\infty)$. 
Thus, (6.28) has no solution for all $n\geq 2$. . Therefore, $\varphi$ is biharmonic if and only if harmonic in this case. 
\par
\underline{Case (3)  $C$-Type}:
 By (6.9), $\varphi:\,M\hookrightarrow {\Bbb C}P^n(4)$ is harmonic if and only if 
$t=\cot u$  $(0<u<\frac{\pi}{4})$ 
must satisfy 
\begin{equation}
(n-2)t^4-2(n+2)t^2+n-2=0,
\end{equation}
which is equivalent to that 
\begin{equation}
t=\cot u=\left\{
\frac{n+2\pm2\sqrt{2n}}{n-2}
\right\}^{1/2}
=\frac{\sqrt{n}\pm\sqrt{2}}{\sqrt{n-2}}.
\end{equation}
On the other hand, by Theorem 5.1 and (6.10), 
$\varphi:\,M\hookrightarrow {\Bbb C}P^n(4)$ is non-harmonic but biharmonic if and only if 
\begin{align}
g(X):&=(n-2)X^2(X-1)^2+(n-2)(X-1)^2+4X(X^2+6X+1)\nonumber\\
&-2X(X-1)^2-2(n+1)X(X-1)^2=0,
\end{align}
where $X:=t^2$. But, 
$g(X)>0$ for all $0<X<\infty$ and $n\geq 3$. Indeed, 
$(1)$ we have 
$$
g(X)=(n-2)(X-1)^2\left\{
X^2-2\frac{n+2}{n-2}X+1
\right\}+4X(X^2+6X+1),
$$
which is positive when either
$X> 5+2\sqrt{6}$ or $0<X<5-2\sqrt{6}$
if $n\geq 3$. 
Furthermore, (2) 
we have 
$$
g(X)=(n-2)(X-1)^4+4X\left(-X^2+10X-1\right),
$$
and 
$-X^2+10X-1>0$ if 
$5-2\sqrt{6}<X<5+2\sqrt{6}$. 
Finally, (3) 
we have 
$g(5\pm 2\sqrt{6})=(4\pm2\sqrt{6})^4>0$. 
Thus, by (1), (2) and (3), $g(X)>0$ on $(0,\infty)$ when $n\geq 3$. 
Thus, (6.31) has no solution for all $n\geq 3$. . Therefore, $\varphi$ is biharmonic if and only if harmonic in this case. 
\par
\underline{Case (4) $D$-type}. 
By (6.13), $\varphi:\,M\hookrightarrow{\Bbb C}P^9$ 
is harmonic if and only if $t=\cot u=\frac{1}{5}$, and 
by (6.14), 
is biharmonic but not harmonic if and only if 
$t=\cot u$ is a solution of the equation 
\begin{equation}
11X^3+63X^2+X+5-20X(X-1)^2=0
\end{equation}
which is equivalent to 
\begin{equation}
h(X):=11X^3+43X^2+41X-15=0.
\end{equation}
This has a solution because 
$h(0)=-15<0$, $h(X)>0$ for a large $X$, and the mean value theorem. 
 Indeed, The solution $X$ of (6.33) is $0.278629$, and 
 the corresponding $t=\cot u$ is $0.527853$, and $u$ is $1.08512$. 
 \par
 \underline{Case (5) $E$-type}. 
By (6.17), $\varphi:\,M\hookrightarrow{\Bbb C}P^{15}$ 
is harmonic if and only if 
$t=\cot u=\frac{\sqrt{15}\pm \sqrt{6}}{3}$
if and only if $u$ is $0.443039$ or $1.12776$.
By (6.18), 
is biharmonic but not harmonic if and only if 
$t=\cot u$ is a solution of the equation 
\begin{equation}
21X^3+99X^2-9X+9-2X(X-1)^2=0
\end{equation}
which is equivalent to 
\begin{equation}
k(X):=13X^3-107X^2+43X-9=0.
\end{equation}
This has a solution because 
$k(0)=-9<0$, $k(X)>0$ for a large $X$, and the mean value theorem. 
 Indeed, The solution $X$ of (6.35) is $7.81906$, and 
 the corresponding $t=\cot u$ is $2.79626$, and $u=0.343448$. 
\end{pf}
\vskip0.6cm\par
\section{Biharmonic homogeneous real hypersurfaces 
in the quarternionic projective space}
In this section, we show classification of all the real hypersurfaces curvature adapted in the quarternionic projective space 
${\Bbb H}P^n(4)$ which are  {\em biharmonic}. 
\par
Let $(N,h)={\Bbb H}P^n(c)$ be the quaternionic projective space with quarternionic sectional curvature $c>0$. Then, the Riemannian curvature tensor is given by 
\begin{align}
R&(U,V)W=\frac{c}{4}
\bigg\{
h(V,W)U-h(U,W)V\nonumber\\
&\qquad
+\sum_{\alpha=1}^3
\big(h(J_{\alpha}V,W)J_{\alpha}U
-h(J_{\alpha}U,W)J_{\alpha}V
+2h(U,J_{\alpha}V)J_{\alpha}W\big)\bigg\},\nonumber
\end{align}
for vector fields $U$, $V$ and $W$ on ${\Bbb H}P^n(c)$. 
Here, $J_{\alpha}$ $(\alpha=1,2,3)$ are 
the locally defined adapted three almost complex tensors on ${\Bbb H}P^n(c)$ 
which satisfy 
$J_1J_2=-J_2J_1=J_3$. 
Then, we have the following 
theorem which we omit its proof since 
one can prove it by the same manner as Theorem 5.1 whose proof is ommited.
\begin{th}
Let $(M,g)$ be a real $(4n-1)$-dimensional compact Riemannian manifold, and $\varphi:\,(M,g)\rightarrow {\Bbb H}P^n(c)$ be an isometric immersion with constant non-zero mean curvature 
$(n\geq 2)$. Then, the necessary and sufficient condition for $\varphi$ to be biharmonic is 
\begin{equation}
\Vert B(\varphi)\Vert^2=(n+2)c.
\end{equation}
\end{th}
\vskip0.6cm\par
Now, let us recall Berndt's classification (\cite{B}) of 
all the real hypersurfaces $(M,g)$ 
in the quarternionic projective space 
${\Bbb H}P^n(4)$ 
which are 
{\em curvature adapted}, i.e., 
$J_{\alpha}\xi$ is a direction of the principal curvature for all 
$\alpha=1,2,3$, where $\xi$ is the unit normal vector field 
along $M$.  
\begin{th} $($Berndt$ \,\,\cite{B})$ 
$({\rm I})$ All the curvature adapted real hypersurfaces in ${\Bbb H}P^n(4)$ are one of the following:
\par
$(1)$ a geodesic sphere $M(u)$ of radius $u$ $(0<u<\frac{\pi}{2})$,
\par
$(2)$ a tube $M(u)$ of radius $u$ $(0<u<\frac{\pi}{4})$ of 
the complex projective space ${\Bbb C}P^n\subset {\Bbb H}P^n(4)$, and 
\par
$(3)$ tubes $M_k(u)$ of radii $u$ $(0<u<\frac{\pi}{4})$ 
of  the quaternionic projective subspaces 
${\Bbb H}P^k\subset {\Bbb H}P^n(4)$ with $1\leq k\leq n-1$. 
\vskip0.3cm \par
 $({\rm II})$ Furthermore, their principal curvatures are given as follows. 
 \par
 $(1)$ The geodesic sphere $M(u)$:
 \begin{equation}
 \left\{
 \begin{aligned}
 \lambda_1&=\cot u \,\,{\rm (with \,\,multiplicity} \,\, m_1=4(n-1)), \\
 \lambda_2&=2\cot (2u)\,\, {\rm (with \,\,multiplicity}\,\, m_2=3). 
 \end{aligned}
 \right.
 \end{equation}
 \par
 $(2)$ The tube $M(u)$ of the complex projective space:
 \begin{equation}
 \left\{
 \begin{aligned}
 \lambda_1&=\cot u \,\,{\rm (with \,\,multiplicity} \,\, m_1=2(n-1)), \\
 \lambda_2&=-\tan u\,\, {\rm (with\,\,multiplicity}\,\, m_2=2(n-1)),\\
 \lambda_3&=2\cot(2u)\,\,{\rm (with \,\,multiplicity}\,\,m_3=1),\\
 \lambda_4&=-2\tan(2u)\,\,{\rm (with \,\,multiplicity}\,\,m_4=2).
 \end{aligned}
 \right.
 \end{equation}
 \par
 $(3)$ The tubes $M_k(u)$ of the quarternionic projective spaces: 
 \begin{equation}
 \left\{
 \begin{aligned}
 \lambda_1&=\cot u \,\,{\rm (with\,\,multiplicity} \,\, m_1=4(n-k-1)), \\
 \lambda_2&=-\tan u\,\, {\rm (with \,\,multiplicity}\,\, m_2=4k),\\
 \lambda_3&=2\cot(2u)\,\,{\rm (with \,\,multiplicity}\,\,m_3=3).\\
  \end{aligned}
 \right.
 \end{equation}
\end{th}
\vskip0.6cm\par
Then, we obtain the following theorem. 
\begin{th}
For all the three classes $(1)$, $(2)$ and $(3)$ of Theorem 7.2, 
harmonic (i.e, minimal), and 
biharmonic but not harmonic real hypersurfaces 
$M(u)$ or $M_k(u)$ in ${\Bbb H}P^n(4)$ with radii $u$ 
are given as follows:
\par
$(1)$ The geodesic sphere $M(u)$$:$ \quad 
The necessary and sufficient condition for $M(u)$ is 
to be harmonic $($i.e., minimal$)$ is that
$t=\cot u$ $(0<u<\frac{\pi}{2})$ satisfies 
\begin{equation}
t=\sqrt{\frac{3}{4n-1}},
\end{equation}
and 
to be biharmonic  but not harmonic is that 
$t=\cot u$ $(0<u<\frac{\pi}{2})$ satisfies 
\begin{equation}
(4n-1)t^4-2(2n+7)t^2+3=0.
\end{equation}
Both the $(7.5)$ and $(7.6)$  have always solutions.
\par
$(2)$ The tube $M(u)$ of radius $u$ $(0<u<\frac{\pi}{4})$ 
of the complex projective space$:$ \quad
The necessary and sufficient condition for $M(u)$ is to be 
harmonic $($i.e., minimal$)$ is that 
\begin{equation}
(2n-1)t^4-(4n+5)t^2+2(n-1)=0, 
\end{equation}
and to be 
biharmonic but not harmonic is that 
\begin{equation}
(2n-1)t^8-8(n+1)t^6-(6n+11)t^4-2(2n-1)t^2-12=0.
\end{equation}
Both the $(7.7)$ and $(7.8)$ have always solutions. 
\par
$(3)$ The tubes $M_k(u)$ of radii $u$ $(0<u<\frac{\pi}{4})$ of the quarternioinic projective subspaces$:$
\quad
The necessary and sufficient conditions for $M_k(u)$ to be 
harmonic $($i.e., minimal$)$ is that 
\begin{equation}
t=\sqrt{\frac{4k+3}{4n-4k-1}},
\end{equation}
and to be biharmonic but not harmonic is that 
\begin{equation}
(4n-4k-1)t^4-2(2n+4)t^2+4k+3=0.
\end{equation}
Both the $(7.9)$ and $(7.10)$ have always solutions. 
\end{th}
\begin{pf}
\underline{Case (1): The geodesic sphere $M(u)$}. 
\quad In this case, the mean curvature $H$ of $M(u)$ is given by 
\begin{align}
H&=\frac{1}{4n-1}\left\{
4(n-1)\cot u+3\,2\cot (2u)
\right\}\nonumber\\
&=4(n-1)t+3\left(
t-\frac{1}{t}
\right),
\end{align}
where $t=\cot u$, so that $M(u)$ is harmonic, i.e., minimal
if and only if 
\begin{equation}
4(n-1)t^2+3t -3=0\quad\Longleftrightarrow \quad 
t=\sqrt{\frac{3}{4n-1}}.
\end{equation}
The square of the second fundamental form 
$\Vert B(\varphi)\Vert ^2$ is given by 
\begin{equation}
\Vert B(\varphi)\Vert^2=
4(n-1)t^2+3\left(
t-\frac{1}{t}
\right)^2
=(4n-1)t^2+\frac{3}{t^2}-6,
\end{equation}
which yields by Theorem 7.1, that 
$M(u)$ is biharmonic, but not harmonic  if and only if 
\begin{align}
(4n-1)&t^2+\frac{3}{t^2}-2(2n+7)=0\nonumber\\
&\Longleftrightarrow\quad 
t^2=\frac{2n+7\pm\sqrt{n^2+4n+13}}{4n-1},
\end{align}
which has always solutions. 
\par
\underline{Case (2): The tube $M(u)$ 
of ${\Bbb C}P^n\subset {\Bbb H}P^n(4)$}. 
\quad 
In this case, the mean curvature $(4n-1)H$ of  $M(u)$ coincides with  
\begin{align}
&2(n-1)\cot u+2(n-1)(-\tan u)+2\cot(2u)+2(-2\cot(2u))\nonumber\\
&=2(n-1)+2(n-1)\left(-\frac{1}{t}\right)
+\left(t-\frac{1}{t}
\right)+2\left(\frac{-4t}{t^2-1}
\right)\nonumber\\
&=\frac{(2n-1)t^4-(4n+5)t^2+2(n-1)}{t(t^2-1)}, 
\end{align}
where $t=\cot u$, so that $M(u)$ is harmonic, i.e., minimal 
if and only if 
\begin{align}
2(n-1)t^4&-(4n+5)t^2+2(n-1)=0\nonumber\\
&\Longleftrightarrow
\quad 
t^2=\frac{4n+5\pm\sqrt{3(n+2)(2n+9)}}{2(n-1)},
\end{align}
which has always solutions. 
On the other hand, $\Vert B(\varphi)\Vert ^2$ coincides with 
\begin{align}
&2(n-1)t^2+\frac{2(n-1)}{t^2}+t^2-2+\frac{1}{t^2}+\frac{32}{(t^2-1)^2}\nonumber\\
&=\frac{(2n-1)X^2(X-1)^2+(2n-1)(X-1)^2-2X(X-1)^2
+32X^2}{X(X-1)^2}
\end{align}
where $X=t^2$. Hence, $M(u)$ is biharmonic, but not harmonic if and only if 
\begin{equation}
(2n-1)X^4-8(n+1)X^3-(6n+11)X^2-2(2n-1)X-12=0, 
\end{equation}
with $X=t^2$, $t=\cot u$ with $0<u<\frac{\pi}{4}$. 
 Denoting by $f(t)$ the LHS, 
 $$
 f(0)=-12<0,\quad f(t)>0
 $$
 for large $t$. Thus, by the mean value theorem, (7.18) has always solutions $X$, so $t$, but not solutions of (7.16). 
\par
\underline{Case (3): The tubes of ${\Bbb H}P^k\subset {\Bbb H}P^n(4)$}. 
\quad 
In this case, the mean curvature $H$ of $M(u)$ is given by 
\begin{align}
H&=
4(n-k-1)\cot u+4k (-\tan u)+6\cot(2u)\nonumber\\
&=
4(n-k-1)t+4k\left(-\frac{1}{t}\right)+
3\left(t-\frac{1}{t}\right),
\end{align}
with $t=\cot u$, so that $M(u)$ is harmonic, i.e., minimal if and only if 
\begin{equation}
t=\sqrt{\frac{4k+3}{4n-4k-1}}.
\end{equation}
On the other hand, $\Vert B(\varphi)\Vert ^2$ is given by 
\begin{align}
\Vert B(\varphi)\Vert^2
&=4(n-k-1)\cot^2 u+4k\tan^2u+12\cot^2(2u)\nonumber\\
&=(4n-4k-1)t^2+\frac{4k}{t^2}+3\left(t-\frac{1}{t}\right)^2\nonumber\\
&=(4n-4k-1)t^2+\frac{4k+3}{t^2}-6,
\end{align}
so that $M(u)$ is biharmonic, but not harmonic if and only if 
\begin{equation}
(4n-4k-1)t^4-2(2n+4)t^2+4k+3=0, 
\end{equation}
which has always solutions. 
\end{pf}
\vskip0.6cm\par
\section{Biharmonic maps into a manifold of nonpositive curvature}
In this section, we show answers in case of bounded geometry, to the following 
conjectures proposed by B.Y. Chen (\cite{C}), and 
R. Caddeo, S. Montaldo and P. Piu (\cite{CMP}): 
\par
{\bf B.Y. Chen's Conjecture.}
{\em Any biharmonic submanifold of the Euclidean space is harmonic.}
\par\noindent or more generally, \par
{\bf R. Caddeo, S. Montaldo and P. Piu's conjecture.}  
{\em The only biharmonic submanifolds of a complete Riemanian manifold whose curvature is nonpositive are the minimal ones.}
\vskip0.6cm\par
{\bf Example 8.1} \quad 
Let $\varphi:\,({\mathbb R}^m,g_0)\ni x=(x_1,\dots,x_m)\mapsto (\varphi_1,\dots,\varphi_n)\in ({\mathbb R}^n,h_0)$ be a smooth mapping given by 
\begin{equation}
\varphi_i(x)=\sum_{j=1}^mx_j{}^4-m\,x_i{}^4\quad (i=1,\cdots,m),\nonumber
\end{equation}
and $\varphi_j(x)$ $(j=m+1,\dots,n)$ are at most linear, 
where $({\mathbb R}^m,g_0)$ and $({\mathbb R}^n,h_0)$ are the standard Euclidean spaces, respectively. 
Then, we have 
\begin{equation}\left\{
\begin{aligned}
\tau(\varphi)&=\Delta\varphi=(\Delta\varphi_1,\dots,\Delta\varphi_n),\\
\tau_2(\varphi)&=\Delta(\Delta\varphi)=0,
\end{aligned}
\right.\nonumber
\end{equation}
where 
\begin{equation}
\Delta\varphi_i=12\left(\sum_{j=1}^mx_j{}^2-m\,x_i{}^2\right)
\quad (i=1,\dots,m).\nonumber
\end{equation}
Furthermore, we have
\begin{align}
\Vert \tau(\varphi)\Vert^2 
&=12^2\,m\left(m\sum_{j=1}^mx_j{}^4-\left(\sum_{j=1}^mx_j{}^2\right)^2\right)\geq 0,\nonumber\\
\Vert\overline{\nabla}\tau(\varphi)\Vert^2
&=24^2\,m(m-1)\left(\sum_{j=1}^mx_j{}^2\right)^2.\nonumber
\end{align}
\vskip0.6cm\par
However, we show
\begin{th} 
Let $\varphi:\,(M,g)\rightarrow (N,h)$ be a biharmonic map 
from a complete Riemannian manifold $(M,g)$ of bounded sectional curvature $\vert{\rm Riem}^M\vert \leq C$ into a Riemannian manifold 
$(N,h)$ of nonpositive curvature, i.e., ${\rm Riem}^N\leq 0$. 
Assume that 
\begin{equation}
\Vert \tau(\varphi)\Vert\in L^2(M), \,\,\text{and}\quad 
\Vert \overline{\nabla}\tau(\varphi)\Vert\in L^2(M).
\end{equation}
Then, $\varphi:\, (M,g)\rightarrow (N,h)$ is harmonic. 
\end{th}
\vskip0.6cm\par
\begin{cor}
Let $\varphi:\, (M,g)\rightarrow (N,h)$ be a biharmonic isometric immersion
from a complete Riemannian manifold $(M,g)$ of bounded sectional curvature $\vert{\rm Riem}^M\vert \leq C$ into a Riemannian manifold 
$(N,h)$ of nonpositive curvature, i.e., 
${\rm Riem}^N\leq 0$.
Assume that the second fundamental form 
$\tau(\varphi)$ satisfies that 
\begin{equation}
\Vert \tau(\varphi)\Vert\in L^2(M), \,\,\text{and}\quad
\Vert \overline{\nabla}\tau(\varphi)\Vert\in L^2(M).
\end{equation}
Then $\varphi:\,(M,g)\rightarrow (N,h)$ is harmonic. 
\end{cor}
\vskip0.6cm\par
Before going to prove Theorem 8.1, we prepare a cut off function 
$\lambda_R$ $(0<R<\infty)$ 
on a complete Riemannian manifold $(M,g)$ as follows (\cite{DMO}). 
Let $\mu$ be a real valued $C^{\infty}$ function on ${\Bbb R}$ 
satistying the following conditions:
\begin{equation}
\left\{
\begin{aligned}
&0\leq \mu(t)\leq 1\quad (t\in {\Bbb R}),\\
&\mu(t)=1\qquad (t\leq 1),\\
&\mu(t)=0\qquad (t\geq 2),\\
&\vert \mu'\vert\leq C,\,\, {\rm and}\,\,
\vert\mu''\vert\leq C,
\end{aligned}
\right.
\end{equation}
where $\mu'(t)$ and $\mu''(t)$ stand for the derivations of the first and second order of $\mu(t)$ with respect to $t$, respectively. Then, 
for all $R>0$, the function defined by 
$$
\lambda_R(x)=\mu\left(
\frac{r(x)}{R}
\right), \quad (x\in M)
$$
is said to be a {\em cut off function} on $(M,g)$, where 
$$
r(x)=d(x_0,x),\quad (x\in M)
$$
for some fixed point $x_0$ in $M$ and $d(x,y)$, $(x,y\in M)$ is the Riemannian distance function of $(M,g)$.  
Then, it is known (\cite{DMO}) that 
\begin{lem}
$({\rm i})$ \quad $\lambda_R$ is a Lipshitz function on $M$, and differentiable a.e. on $M$, 
\par
$({\rm ii})$ \quad ${\rm supp}(\lambda_R)\subset B_{2R}(x_0)$, 
\par
$({\rm iii})$ \quad $0\leq \lambda_R(x)\leq 1$, \quad $(x\in M)$, 
\par
$({\rm iv})$ \quad $\lambda_R(x)=1$, \quad $(x\in B_R(x_0))$, 
\par
$({\rm v})$ \quad $\vert\nabla\lambda_R\vert\leq \frac{C}{R}$, $($a.e. on 
$M)$, 
\par
$({\rm vi})$ \quad and 
if the Ricci curvature of $(M,g)$ is bounded below by 
a constant $(m-1)(-k)$ for some $k>0$ $(m=\dim M)$, then, 
\begin{equation}
\vert \Delta \lambda_R\vert \leq
\frac{C}{R^2}+\frac{CC'}{R}\quad ({a.e.\,\, on}\,\, M).
\end{equation}
Here $C'$ is a positive constant depending only on $m$ and $k$, 
${\rm supp}(\lambda_R)$ stands for the support of $\lambda_R$, 
and $B_r(x):=\{y\in M;\,d(x_,y)<r\}$ 
is the Riemannian disc in $(M,g)$ around $x$ with radius $r>0$. 
\end{lem}
\begin{pf}
From $({\rm i})$ to $({\rm v})$, see \cite{DMO}, for instance. 
For $({\rm vi})$, let us recall the estimation of $\Delta r$ in terms of the lower bound of the Ricci curvature (see \cite{Ka} for instance): 
If the Ricci curvature of $(M,g)$ is bounded below by 
a constant $(m-1)(-k)$ for some $k>0$ $(m=\dim M)$, then, 
\begin{equation}
\Delta r
\leq 
(m-1)\frac{f_k{}'}{f_k}
=\frac{m-1}{\sqrt{k}}\,
\frac{\cosh(\sqrt{k}r)}{\sinh(\sqrt{k} r)},
\end{equation}
where where $f_k(t)=\frac{\sinh (\sqrt{k}t)}{\sqrt{k}}$ is the unique solution of the initial value problem 
$$
f_k{}''+(-k)f_k=0,\quad f_k(0)=0,\,\,f_k'(0)=1. 
$$
Thus, outside of $B_R(x_0)$, it holds that 
\begin{equation}
\vert\Delta r\vert \leq 
\frac{m-1}{\sqrt{k}}\,
\frac{\cosh(\sqrt{k}R)}{\sinh(\sqrt{k}R)}.
\end{equation}
Since 
$\nabla \lambda_R=\frac{1}{R}\,\mu'\left(\frac{r}{R}
\right)\,\nabla r$ 
(see \cite{Ka}, p. 108), 
we have, a.e. on $M$, 
\begin{equation}
\Delta \lambda_R=\frac{1}{R^2}\,\mu''\left(\frac{r}{R}\right)
+\frac{1}{R}\,\mu'\left(\frac{r}{R}\right)\,\Delta r.
\end{equation}
Then, together with (8.3), (8.6) and (8.7), we have (8.4). 
\end{pf}
\vskip0.6cm\par
Now let us begin a proof of Theorem 8.1. 
Let us recall the definition of 
$e_2(\varphi)=\frac12\Vert \tau(\varphi)\Vert ^2$. We will estimate 
$\Delta(\lambda_R\,e_2(\varphi))$ as follows. 
 \begin{equation}
 \Delta (\lambda_R \,e_2(\varphi))
 =(\Delta \lambda_R)\,e_2(\varphi)+2g(\nabla \lambda_R,\nabla e_2(\varphi))+\lambda_R\, \Delta e_2(\varphi).
 \end{equation}
 \par
 For the LHS of (8.8), 
 we have
 $\Delta(\lambda_R \,e_2(\varphi))={\rm div}X$, 
 where 
 $X:=\nabla(\lambda_R\,e_2(\varphi))$ which  
 is a $C^{\infty}$ vector field on $M$ with compact support.   
 Due to Green's theorem, 
 \begin{equation}
 \int_M\Delta(\lambda_R\,e_2(\varphi))v_g=
 \int_M{\rm div}(X)v_g=0. 
 \end{equation}
 Furthermore, we have 
 \begin{align}
 \lim_{R\rightarrow \infty}\int_M(\Delta\lambda_R)e_2(\varphi)v_g=0,\\
 \lim_{R\rightarrow \infty}\int_Mg(\nabla\lambda_R,\nabla e_2(\varphi))v_g=0.
 \end{align}
 Indeed, for (8.10), by (8.4) in Lemma 8.1, 
 \begin{align}
 \bigg\vert
 \int_M(\Delta\lambda_R)e_2(\varphi) v_g\bigg\vert
 &\leq 
 \int_M\vert\Delta\lambda_R\vert\,e_2(\varphi)v_g\nonumber\\
 &\leq \int_M\left(\frac{C}{R^2}+\frac{CC'}{R}\right)\,e_2(\varphi)v_g\nonumber\\
 &=\left(\frac{C}{R^2}+\frac{CC'}{R}\right)\,\int_Me_2(\varphi)v_g,
 \end{align}
 where the RHS goes to $0$ if $R\rightarrow \infty$, since 
 $e_2(\varphi)=\frac12\Vert \tau(\varphi)\Vert^2\in L^1(M)$ by the assumptions (8.1). 
 For (8.11), due to $({\rm v})$ in Lemma 8.1, 
 \begin{align}
 \bigg\vert \int_Mg(\nabla\lambda_R,\nabla e_2(\varphi))v_g\bigg\vert
 &\leq \int_M\vert g(\nabla\lambda_R,\nabla e_2(\varphi))\vert v_g\nonumber\\
 &\leq 
 \int_M\Vert \nabla\lambda_R\Vert\,\,\Vert \nabla e_2(\varphi)\Vert v_g\nonumber\\
 &\leq \frac{C}{R}\int_M\Vert \nabla e_2(\varphi)\Vert v_g,
 \end{align}
 where the RHS goes to $0$ if $R\rightarrow \infty$, since 
 \begin{align}
 \int_M\Vert \nabla e_2(\varphi)\Vert v_g&=\int_M\vert g(\overline{\nabla}\tau(\varphi), \tau(\varphi))\vert v_g\nonumber\\
 &\leq\int_M\Vert \overline{\nabla}\tau(\varphi)\Vert\,\Vert\tau(\varphi)\Vert v_g\nonumber\\
 &\leq\Vert \overline{\nabla}\tau(\varphi)\Vert_{L^2(M)}
 \Vert\tau(\varphi)\Vert_{L^2(M)}<\infty\nonumber
 \end{align}
 by the assumptions (8.1). 
 \par
 Thus, due to ((8.8), (8.9), (8.10), (8.11), we obtain
 \begin{equation}
 \lim_{R\rightarrow\infty}\int_M\lambda_R\,\Delta e_2(\varphi)v_g=0.
 \end{equation}
 \par
 Now, by the computation (4.1) in \cite{J} in which Jiang used only 
 the assumption that $\varphi:\,(M,g)\rightarrow (N,h)$ is biharmonic, 
 we have 
 \begin{align}
 \Delta e_2(\varphi)&=
\sum_{k=1}^m
g(\overline{\nabla}_{e_k}\tau(\varphi),\overline{\nabla}_{e_k}\tau(\varphi))
+g(-\overline{\nabla}^{\ast}\overline{\nabla}\tau(\varphi),\tau(\varphi))\nonumber\\
&=
\sum_{k=1}^m
g(\overline{\nabla}_{e_k}\tau(\varphi),\overline{\nabla}_{e_k}\tau(\varphi))\nonumber\\
&\,\,-\sum_{k=1}^mh(R^N(\tau(\varphi),d\varphi(e_k))d\varphi(e_k),\tau(\varphi)).
 \end{align}
 Then, we have 
  \begin{align}
 \int_M&\lambda_R\,\Delta e_2(\varphi)v_g
 =\int_M\lambda_R\left(
 \sum_{k=1}^m
g(\overline{\nabla}_{e_k}\tau(\varphi),\overline{\nabla}_{e_k}\tau(\varphi))\right)
v_g\nonumber\\
&\,\,+\int_M\lambda_R
\left(-\sum_{k=1}^mh(R^N(\tau(\varphi),d\varphi(e_k))d\varphi(e_k),\tau(\varphi))\right)v_g. 
\end{align}
Here, the first term of the RHS of (8.16) goes to $0$ when $R\rightarrow\infty$, i.e.,  
\begin{equation}
\lim_{R\rightarrow\infty}\int_M\lambda_R\left(
 \sum_{k=1}^m
g(\overline{\nabla}_{e_k}\tau(\varphi),\overline{\nabla}_{e_k}\tau(\varphi))\right)
v_g=0.
\end{equation}
Because, both the integrand of (8.17) is nonnegative, and  
by the curvature assumption of $(N,h)$, ${\rm Riem}^N\leq 0$, 
the integrand of the second term of RHS of (8.16) is nonnegative. 
Thus, (8.14) implies the desired (8.17). 
\par
Notice here, that (8.17) implies 
\begin{equation}
\int_M\sum_{k=1}^mg(\overline{\nabla}_{e_k}\tau(\varphi),
\overline{\nabla}_{e_k}\tau(\varphi))v_g=0,
\end{equation}
which yields that 
$
\overline{\nabla}_X\tau(\varphi)=0
$
for all $X\in {\frak X}(M)$. 
\par
Finally, if we consider a $C^{\infty}$ vector field $X_{\varphi}$ on $M$ defined by
$$
X_{\varphi}:=\sum_{k=1}^mh(d\varphi(e_k),\tau(\varphi))e_k,
$$
the divergence of $X_{\varphi}$ satisfies that 
\begin{align}
{\rm div}(X_{\varphi})&=
h(\tau(\varphi),\tau(\varphi))+\sum_{k=1}^mh(d\varphi(e_k),\overline{\nabla}_{e_k}\tau(\varphi))\nonumber\\
&=h(\tau(\varphi),\tau(\varphi))\in L^1(M),
\end{align}
by the above and the assumptions (8.1). Therefore, due to the Green's theorem on a complete Riemannian manifolds $(M,g)$ (see \cite{G} for instance), we obtain 
\begin{equation}
\int_Mh(\tau(\varphi),\tau(\varphi))v_g=
\int_M{\rm div}(X_{\varphi})v_g=0,
\end{equation}
which yields $\tau(\varphi)=0$. 
\qed
\vskip0.6cm\par
\section{The first variational formula for bi-Yang-Mills fields}
From this section, we begin to prepare fundamental materials to state 
interesting phenomena on  bi-Yang-Mills fields 
which are closely related to biharmonic maps. 
We will recall the Yang-Mills setting (\cite{BL}) and the definition of bi-Yang-Mills fields following Bejan and Urakawa (\cite{BU}), and show the isolation phenomena. 
\vskip0.3cm\par
Let us start with the Yang-Mills setting following \cite{BL}. Let $(E,h)$ be a real vector bundle of rank $r$ with an inner product $h$ 
over an $m$-dimensional compact Riemannian manifold $(M,g)$. 
Let ${\mathcal C}(E,h)$ be the space of all $C^{\infty}$-connections of $E$ satisfying the compatibility condition: 
$$
Xh(s,t)=h(\nabla_Xs,t)+h(s,\nabla_Xt),\quad s,\,t\in \Gamma(E),
$$
for all $X\in {\frak X}(M)$, where $\Gamma(E)$ stands for the space of all $C^{\infty}$-sections of $E$. For $\nabla\in {\mathcal C}(E,h)$, let $R^{\nabla}$ be its curvature tensor defined by 
$$
R^{\nabla}(X,Y)s=\nabla_X(\nabla_Ys)-\nabla_Y(\nabla_Xs)-\nabla_{[X,Y]}s, 
$$
for all $X,Y\in {\frak X}(M),\,s\in \Gamma(E)$. Let $F={\rm End}(E,h)$ be the bundle of endmorphisms of $E$ which are skew symmetric with respect to the inner product $h$ on $E$. We define the inner product 
$\langle \,,\,\rangle$ on $F$ by 
$$
\langle \varphi,\psi\rangle=\sum_{i=1}^rh(\varphi u_i,\psi u_i), \quad\varphi,\psi\in F_x,
$$
where $\{u_i\}_{i=1}^r$ is an orthonormal basis of $E_x$ with respect to $h$ $(x\in M)$. Let us also consider the space of $F$-valued $k$-forms on $M$, denoted by $\Omega^k(F)=\Gamma(\wedge^kT^{\ast}M)\otimes F)$, which admits a global inner product $(\,,\,)$ given by 
$$
(\alpha,\beta)=\int_M\langle\alpha,\beta\rangle v_g,
$$
where the pointwise inner product $\langle\alpha,\beta\rangle$ is given by 
$$
\langle\alpha,\beta\rangle=\sum_{i_1<\cdots<i_k}
\langle\alpha(e_{i_1},\dots,e_{i_k}),\beta(e_{i_1},\dots,e_{i_k})\rangle
$$
and $\{e_i\}_{i=1}^m$ is a locally defined orthonormal frame field on 
$(M,g)$. 
\par
For every $\nabla\in {\mathcal C}(E,h)$, let $d^{\nabla}:\,\Omega^k(F)\rightarrow\Omega^{k+1}(F)$ be the exterior differentiation with respect to $\nabla$ (cf. \cite{BL}), and the adjoint operator 
$\delta^{\nabla}:\,\Omega^{k+1}(F)\rightarrow \Omega^{k}(F)$ given by 
$$
\delta^{\nabla}\alpha=(-1)^{k+1}\ast d^{\nabla}\ast \alpha,\quad\alpha\in \Omega^{k+1}F),
$$
where $\ast:\,\Omega^p(F)\rightarrow\Omega^{m-p}(F)$ is the extension of the usual Hodge star operator on $(M,g)$. Then, it holds that 
$$
(d^{\nabla}\alpha,\beta)=(\alpha,\delta^{\nabla}\beta),\quad \alpha\in \Omega^k(F),\,\beta\in \Omega^{k+1}(F). 
$$
\vskip0.3cm\par
Now let us recall the bi-Yang-Mills functional (see \cite{BU}) and Yang-Mills one (see \cite{BL}):
\begin{df}
\begin{align}
{\mathcal YM}_2(\nabla)&=\frac12\int_M\Vert\delta^{\nabla}R^{\nabla}\Vert^2v_g,\quad \nabla\in {\mathcal C}(E,h), \\
{\mathcal YM}(\nabla)&=\frac12\int_M\Vert R^{\nabla}\Vert^2v_g,\quad \nabla\in {\mathcal C}(E,h), 
\end{align}
where 
$\Vert \delta^{\nabla}R^{\nabla}\Vert$, $($resp. $\Vert R^{\nabla}\Vert$$)$ 
is the norm of $\delta^{\nabla}R^{\nabla}\in \Omega^1(F)$ 
$($resp. $R^{\nabla}\in \Omega^2(F)$$)$ relative to each 
$\langle\,,\,\rangle$. 
\end{df}
\vskip0.6cm\par
Then, the bi-Yang-Mills fields and the Yang-Mills ones are critical points of the above functionals as follows.
\begin{df}
For each $\nabla\in {\mathcal C}(E,h)$, 
it is a {\em bi-Yang-Mills field} $($resp. {\em Yang-Mills field}$)$ 
if for any smooth one-parameter family $\nabla^t$ 
$(\vert t\vert<\epsilon)$ with $\nabla^0=\nabla$, 
\begin{equation}
\frac{d}{dt}\bigg\vert_{t=0}{\mathcal YM}_2(\nabla^t)=0,
\quad\left({\it resp.}\,\,\frac{d}{dt}\bigg\vert_{t=0}{\mathcal YM}(\nabla^t)=0\right).
\end{equation}
\end{df}
\vskip0.6cm\par
Then, the first variation formulas are given as 
\begin{th} $($\cite{BU}, \cite{BL}$)$
Let $\alpha=\frac{d}{dt}\vert_{t=0}\nabla^t\in \Omega^1(F)$. Then, we have 
\begin{align}
\frac{d}{dt}\bigg\vert_{t=0}{\mathcal YM}_2(\nabla^t)&=
\int_M\langle (\delta^{\nabla}d^{\nabla}+{\mathcal R}^{\nabla})
(\delta^{\nabla}R^{\nabla}),\alpha\rangle v_g,\\
\frac{d}{dt}\bigg\vert_{t=0}{\mathcal YM}(\nabla^t)&=
\int_M\langle \delta^{\nabla}R^{\nabla},\alpha\rangle v_g,
\end{align}
respectively. Here, 
${\mathcal R}^{\nabla}(\beta)\in \Omega^1(F)$ $(\beta\in \Omega^1(F))$ is defined by 
\begin{equation}
{\mathcal R}^{\nabla}(\beta)(X)=\sum_{j=1}^m[R^{\nabla}(e_j,X),\beta(e_j)],\quad X\in {\frak X}(M). 
\end{equation}
Thus, $\nabla$ is a bi-Yang-Mills field $($resp. Yang-Mills one$)$ if and only if 
\begin{equation}
(\delta^{\nabla}d^{\nabla}+{\mathcal R}^{\nabla})
(\delta^{\nabla}R^{\nabla})=0\quad ({\it resp.}\,\,
\delta^{\nabla}R^{\nabla}=0).
\end{equation}
\end{th}
\vskip0.6cm\par
Thus, by this theorem, we have immediately 
\begin{cor}
If $\nabla$ is a Yang-Mills field, then it is also a bi-Yang-Mills one. 
\end{cor}
\vskip0.6cm\par
\begin{lem}
For all $\beta_1$, $\beta_2\in \Omega^1(F)$, and $\varphi\in \Omega^2(F)$, we have 
\begin{equation}
\langle \varphi,[\beta_1\wedge\beta_2]\rangle
=\langle {\mathcal R}(\varphi)(\beta_2),\beta_1\rangle
=\langle \beta_2,{\mathcal R}(\varphi)(\beta_1)\rangle.
\end{equation}
\end{lem}
\begin{pf}
For the first equality, we have
\begin{align}
\langle \varphi,[\beta_1\wedge\beta_2]\rangle
&=\sum_{i<j}
\langle 
\varphi(e_i,e_j),[\beta_1\wedge\beta_2](e_i,e_j)\rangle\nonumber\\
&=\sum_{i<j}\varphi(e_i,e_j),[\beta_1(e_i),\beta_2(e_j)]
-[\beta_1((e_j),\beta_2(e_i)]\rangle\nonumber\\
&=\sum_{i,j=1}^m\langle \varphi(e_i,e_j),
[\beta_1(e_i),\beta_2(e_j)]\rangle\nonumber\\
&=\sum_{i=1}^m\left\langle
\sum_{j=1}^m[\varphi(e_j,e_i),\beta_2(e_j)],\beta_1(e_i)
\right\rangle\nonumber\\
&=\sum_{i=1}^m\langle{\mathcal R}(\varphi)(\beta_2)(e_i),\beta_1(e_i)\rangle\nonumber\\
&=\langle{\mathcal R}(\varphi)(\beta_2),\beta_1\rangle, 
\nonumber
\end{align}
since $\langle[\eta,\psi],\xi\rangle+\langle\psi,[\eta,\xi]\rangle=0$ for all endomorphisms 
$\eta$, $\psi$, and $\xi$ of $E_x$ $(x\in M)$. 
By the same reason, for the second equality, we have 
\begin{align}
\langle{\mathcal R}(\varphi)(\beta_2),\beta_1\rangle
&=\sum_{i=1}^m\left\langle
\sum_{j=1}^m[\varphi(e_j,e_i),\beta_2(e_j)],\beta_1(e_i)
\right\rangle\nonumber\\
&=-\sum_{i,j=1}^m
\langle
\beta_2(e_j),[\varphi(e_j,e_i),\beta_1(e_i)]\rangle
\nonumber\\
&=\sum_{j=1}\langle
\beta_2(e_j),\sum_{i=1}^m[\varphi(e_i,e_j),\beta_1(e_i)]\rangle\nonumber\\
&=\sum_{j=1}^m\langle\beta_2(e_j),{\mathcal R}(\varphi)(\beta_1)(e_j)\rangle\nonumber\\
&=\langle \beta_2,{\mathcal R}(\varphi)(\beta_1)\rangle,
\nonumber
\end{align}
thus, we obtain (9.8). 
\end{pf}
\vskip0.6cm\par
\section{Isolation phenomena for bi-Yang-Mills fields}
In this section, we finally show very interesting phenomena which assert 
that Yang-Mills fields are isolated 
among the space of all bi-Yang-Mills fields over compact Riemannian manifolds with positive Ricci curvature.
\begin{th}
$($bounded isolation phenomena$)$ 
Let $(M,g)$ a compact Riemannian 
of which Ricci curvature 
is bounded below by a positive constant 
$k>0$, i.e.,  ${\rm Ric}\geq k\,{\rm Id}$. Assume that 
$\nabla\in {\mathcal C}(E,h)$ is a bi-Yang-Mills field with 
$\Vert R^{\nabla}\Vert< \frac{k}{2}$ pointwisely everywhere on $M$. Then, $\nabla$ is a Yang-Mills field. 
\end{th}
\vskip0.6cm\par
\begin{th} 
$($$L^2$-isolation phenomena$)$ 
Let $(M,g)$ be a four dimensional compact Riemannian manifold of which Ricci curvature is bounded below by a positive constant $k>0$, i.e., ${\rm Ric}\geq k\,{\rm Id}$. 
Assume that ${\nabla}\in {\mathcal C}(E,h)$ is a bi-Yang-Mills field 
satisfying that 
\begin{equation}
\Vert R^{\nabla}\Vert_{L^2}< \frac12\min
\left\{
\frac{\sqrt{c_1}}{18}, \,\,\frac{k}{2}\,{\rm Vol}(M,g)^{1/2}
\right\}.
\end{equation}
Then, $\nabla$ is a Yang-Mills field. 
Here,  $c_1$ is the isoperimetric 
constant of $(M,g)$ given by 
\begin{equation}
c_1=\inf_{W\subset M}\frac{{\rm Vol}_3(W)^4}{
\left(\min\{
{\rm Vol} (M_1),\,{\rm Vol}(M_2)\}\right)^3},
\end{equation}
where $W\subset M$ runs over all the hypersurfaces in $M$, and 
${\rm Vol}_3(W)$ is the three dimensional volume of $W$ with respect to the Riemannian metric on $W$ 
induced from $g$, and 
the complement of $W$ in $M$ has a disjoint union of $M_1$ and $M_2$. 
\end{th}
\vskip0.6cm\par
To prove Theorem 10.1, we need the following Weitzenb\"ock formula.
\begin{lem}
Assume that $\nabla\in {\mathcal C}(E,h)$ is a bi-Yang-Mills field. 
Then, 
\begin{align}
\frac12\Delta\Vert\delta^{\nabla} R^{\nabla}\Vert^2
&=\langle
2{\mathcal R}^{\nabla}(\delta^{\nabla}R^{\nabla})
+\delta^{\nabla}R^{\nabla}\circ{\rm Ric},\delta^{\nabla}R^{\nabla}\rangle\nonumber\\
&\quad+\sum_{i=1}^m\Vert
\nabla_{e_i}(\delta^{\nabla}R^{\nabla})\Vert^2.
\end{align}
Here, 
$\Delta f=\sum_{i=1}^m(e_i{}^2-\nabla_{e_i}e_i)f$ 
is the Laplacian acting on  
smooth functions $f$ on $M$, and, 
for all $\alpha\in \Omega^1(F)$, 
\begin{equation}
(\alpha\circ{\rm Ric})(X):=\alpha({\rm Ric}(X)),\,\,X\in {\frak X}(M),
\end{equation}
where ${\rm Ric}$ is the Ricci transform of $(M,g)$. 
\end{lem}
\begin{pf}
Indeed, for the LHS of (10.3), we have  
\begin{equation}
\frac12\Delta\Vert\delta^{\nabla} R^{\nabla}\Vert^2=\langle-\nabla^{\ast}\nabla(\delta^{\nabla}R^{\nabla}), \delta^{\nabla}R^{\nabla}\rangle+
\sum_{i=1}^m\langle\nabla_{e_i}(\delta^{\nabla}R^{\nabla}),
\nabla_{e_i}(\delta^{\nabla}R^{\nabla})\rangle.
\end{equation}
Let us recall the Weitzenb\"ock formula (cf. \cite{BL}, p.199, Theorem (3.2)) that 
\begin{align}
\Delta^{\nabla}\alpha&=(d^{\nabla}\delta^{\nabla}
+\delta^{\nabla}d^{\nabla})\alpha\nonumber\\
&=\nabla^{\ast}\nabla\alpha+\alpha\circ {\rm Ric}+{\mathcal R}^{\nabla}(\alpha),\quad \alpha\in \Omega^1(F).
\end{align}
It holds that 
\begin{equation}
\delta^{\nabla}(\delta^{\nabla}R^{\nabla})=0. 
\end{equation}
Because for all $\varphi\in \Gamma(F)$, 
$$
(\delta^{\nabla}(\delta^{\nabla}R^{\nabla}),\varphi)=
\int_M\langle R^{\nabla},d^{\nabla}(d^{\nabla}\varphi)\rangle v_g.
$$
But, by using the formula (2.9) in \cite{BL}, p. 194, the integrand 
of the RHS coincides with 
\begin{align}
\langle R^{\nabla},&d^{\nabla}(d^{\nabla}\varphi)\rangle
=\sum_{i<j}\sum_{s=1}^r\langle R^{\nabla}(e_i,e_j)u_s,
(R^{\nabla}(e_i,e_j)\varphi)(u_s)\rangle\nonumber\\
&=
\sum_{i<j}\sum_{s=1}^r\langle 
R^{\nabla}(e_i,e_j)u_s,
(R^{\nabla}(e_i,e_j)(\varphi(u_s))
-\varphi(R^{\nabla}(e_i,e_j) u_s)
\rangle\nonumber\\
&=\sum_{i<j}\langle R^{\nabla}(e_i,e_j),[R(e_i,e_j),\varphi]
\rangle\nonumber\\
&=-\sum_{i<j}\langle[R^{\nabla}(e_i,e_j),R^{\nabla}(e_i,e_j)],\varphi\rangle=
0.\nonumber
\end{align}
since $\langle\psi,[\eta,\xi]\rangle=-\langle[\eta,\psi],\xi\rangle$ for all $\eta,\,\psi,\,\xi\in F={\rm End}(E,h)$.  
\par
Now $\nabla$ is a bi-Yang-Mills field, $(\delta^{\nabla}d^{\nabla}+{\mathcal R}^{\nabla})(\delta^{\nabla}R^{\nabla})=0$, so that 
we have 
\begin{align}
-{\mathcal R}^{\nabla}(\delta^{\nabla}R^{\nabla})
&=\delta^{\nabla}d^{\nabla}(\delta^{\nabla} R^{\nabla})\nonumber\\
&=\Delta^{\nabla}(\delta^{\nabla}R^{\nabla})\quad ({\rm by}\,\,(10.7))\nonumber\\
&=\nabla^{\ast}\nabla(\delta^{\nabla}R^{\nabla})
+\delta^{\nabla}R^{\nabla}\circ {\rm Ric}+{\mathcal R}^{\nabla}(\delta^{\nabla}R^{\nabla}),\nonumber
\end{align}
by (10.6). Thus, we have 
\begin{equation}
-\nabla^{\ast}\nabla(\delta^{\nabla}R^{\nabla})
=2{\mathcal R}^{\nabla}(\delta^{\nabla}R^{\nabla})+\delta^{\nabla}R^{\nabla}\circ {\rm Ric}.
\end{equation}
Substituting (10.8) into the first term of the RHS of (10.5), we have (10.3). 
\end{pf} 
\vskip0.6cm\par
{\it Proof of Theorem 10.1.} \quad
By Integrating (10.3) over $M$, and by Green's theorem, we have
\begin{align}
2\int_M\langle&{\mathcal R}^{\nabla}(\delta^{\nabla}R^{\nabla}),\delta^{\nabla}R^{\nabla}\rangle v_g+\int_M\langle\delta^{\nabla}R^{\nabla}\circ {\rm Ric},\delta^{\nabla}R^{\nabla}\rangle v_g\nonumber\\
&+\int_M\sum_{i=1}^m\langle
\nabla_{e_i}(\delta^{\nabla}R^{\nabla}),
\nabla_{e_i}(\delta^{\nabla}R^{\nabla})\rangle v_g=0.
\end{align}
Notice here that 
\begin{equation}
\vert\langle {\mathcal R}^{\nabla}(\alpha),\alpha\rangle
\vert
\leq \Vert 
R^{\nabla}\Vert\,\Vert\alpha\Vert^2,\quad \alpha\in \Omega^1(F).
\end{equation}
Indeed, by Lemma 9.1, and Schwarz inequality, we have
\begin{align}
\vert\langle{\mathcal R}^{\nabla}(\alpha),\alpha\rangle\vert
&=\vert\langle R^{\nabla},[\alpha\wedge\alpha]\rangle\vert\nonumber\\
&\leq
\Vert R^{\nabla}\Vert\,\Vert[\alpha\wedge\alpha]\Vert\nonumber\\
&\leq \Vert R^{\nabla}\Vert\,\Vert\alpha\Vert^2,
\end{align}
of which the last ienquality follows from 
\begin{align}
\Vert[\alpha&\wedge\alpha]\Vert^2=
\sum_{i<j}\Vert[\alpha\wedge\alpha](e_i,e_j)\Vert^2\nonumber\\
&=\frac12\sum_{i,j=1}^m\Vert[\alpha\wedge\alpha](e_i,e_j)\Vert^2\nonumber\\
&\leq\frac12\sum_{i,j=1}^m
2\Vert\alpha(e_i)\Vert^2\,\Vert\alpha(e_j)\Vert^2
\quad ({\rm Lemma}\,\, (2.30) \,\,{\rm in} \,\,\cite{BL},\,{\rm p}.197)\nonumber\\
&=\left(\sum_{i=1}^m\Vert\alpha(e_i)\Vert^2\right)^2\nonumber
\end{align}
which is equal to $\Vert\alpha\Vert^4$. We have (10.11). 
\par
Furthermore, by the assumption of the Ricci curvature of $(M,g)$, we have 
\begin{equation}
\langle
\delta^{\nabla}R^{\nabla}\circ {\rm Ric},\delta^{\nabla}R^{\nabla}\rangle 
\geq k\Vert\delta^{\nabla}R^{\nabla}\Vert^2.
\end{equation}
Indeed, since, at each point $x\in M$, we may choose an orthonormal basis 
$\{e_i\}_{i=1}^m$ of $(T_xM,g_x)$ in such a way that 
$$
{\rm Ric }(e_i)=\mu_i e_i\quad (i=1\dots,m)
$$
where $\mu_i$ $(i=1,\dots,m)$ are bigger than or equal to $k>0$. 
Then,  
\begin{align}
\langle
\delta^{\nabla}R^{\nabla}\circ {\rm Ric},\delta^{\nabla}R^{\nabla}\rangle &=
\sum_{i=1}^m
\langle
\delta^{\nabla}R^{\nabla}({\rm Ric}(e_i)),\delta^{\nabla}R^{\nabla}(e_i)\rangle\nonumber\\
&=\sum_{i=1}^m\mu_i\Vert\delta^{\nabla}R^{\nabla}(e_i)\Vert^2\nonumber\\
&\geq k\Vert \delta^{\nabla}R^{\nabla}\Vert^2.\nonumber
\end{align}
\par
Under the assumption that 
$\Vert R^{\nabla}\Vert<\frac{k}{2}$ at each point of $M$,  we have 
\begin{equation}
\langle
2{\mathcal R}^{\nabla}(\delta^{\nabla}R^{\nabla})
+\delta^{\nabla}R^{\nabla}\circ{\rm Ric},\delta^{\nabla}R^{\nabla}\rangle\geq 0,
\end{equation}
equality holds if and only if $\delta^{\nabla}R^{\nabla}=0$. 
\par\noindent
Because, by (10.10) and (10.12), we have 
\begin{align}
\langle
2{\mathcal R}^{\nabla}(\delta^{\nabla}R^{\nabla})
&+\delta^{\nabla}R^{\nabla}\circ{\rm Ric},\delta^{\nabla}R^{\nabla}\rangle \nonumber\\
&\geq
(-2\Vert R^{\nabla}\Vert+k)\,\Vert\delta^{\nabla}R^{\nabla}\Vert^2\nonumber\\
&\geq 0,\nonumber
\end{align}
and equality holds if and only if $\Vert\delta^{\nabla}R^{\nabla}\Vert=0$ 
by the assumption 
$\Vert \delta^{\nabla}R^{\nabla}\Vert<\frac{k}{2}$. 
\par
Now due to (10.13), 
both the sum of the first and second terms of th LHS of (10.9), and the third term in the same one are bigger than or equal to $0$. Thus, (10.9) implies that 
the sum of the first and second term of (10.9) is $0$, 
and by (10.13), we have $\delta^{\nabla}R^{\nabla}=0$ 
everywhere on $M$. 
\qed
\vskip0.6cm\par
{\bf Remark 10.1.} 
\quad (1) In the case $\Vert R^{\nabla}\Vert=\frac{k}{2}$, 
we can also conclude 
$\nabla_X(\delta^{\nabla}R^{\nabla})=0$ for all $X\in {\frak X}(M)$. 
\quad
(2) In tha case of the unit sphere $(M,g)=(S^m,{\rm can})$, 
$k=m-1$. 
\vskip0.6cm\par
{\it Proof of Theorem 10.2.} \quad 
For a bi-Yang-Mills field $\nabla\in {\mathcal C}(E,h)$, we have (10.9) which we can estimated by 
(10.10) and (10.12) as follows. 
\begin{align}
0&=2\int_M\langle {\mathcal R}^{\nabla}(\delta^{\nabla}R^{\nabla}), \delta^{\nabla}R^{\nabla}\rangle v_g+\int_M\langle\delta^{\nabla}R^{\nabla}\circ {\rm Ric}, \delta^{\nabla}R^{\nabla}\rangle v_g\nonumber\\
&\qquad+\int_M\Vert\nabla(\delta^{\nabla}R^{\nabla})\Vert^2 v_g\nonumber\\
&\geq \int_M\Vert\nabla(\delta^{\nabla}R^{\nabla})\Vert^2 v_g+k\int_M\Vert\delta^{\nabla}R^{\nabla}\Vert^2 v_g
-2\int_M\Vert R^{\nabla}\Vert\,\Vert \delta^{\nabla}R^{\nabla}\Vert^2 v_g\nonumber\\
&\geq \Vert\nabla(\delta^{\nabla}R^{\nabla})\Vert_{L^2}{}^2
+k\Vert \delta^{\nabla}R^{\nabla}\Vert_{L^2}{}^2-2\Vert R^{\nabla}\Vert_{L^2}\,\Vert\delta^{\nabla}R^{\nabla}\Vert_{L^4}{}^2
\end{align}
by Schwarz inequality. 
\par
Now let us recall (cf. \cite{MO}, p. 160) the Sobolev inequality for a four dimensional Riemannian manifold $(M,g)$:
\begin{equation}
\Vert\nabla f\Vert_{L^2}{}^2\geq \frac{\sqrt{c_1}}{18}\Vert f\Vert_{L^4}{}^2-\frac19\left(\frac{c_1}{{\rm Vol}(M,g)}
\right)^{1/2}\Vert f\Vert_{L^2}{}^2,\quad f\in H^2_1(M),
\end{equation}
where $H^2_1(M)$ is the Sobolev space of $(M,g)$. 
By applying (10.15) to the first term of (10.14), 
we have 
\begin{align}
{\rm the\,\, RHS\,\,of\,\,(10.14)}&\geq 
\frac{\sqrt{c_1}}{18}\Vert \delta^{\nabla}R^{\nabla}\Vert_{L^4}{}^2
-\frac19\left(\frac{c_1}{{\rm Vol}(M,g)}\right)^{1/2}\Vert \delta^{\nabla}R^{\nabla}\Vert_{L^2}{}^2\nonumber\\
&\qquad+k\Vert \delta^{\nabla}R^{\nabla}\Vert_{L^2}{}^2
-2\Vert R^{\nabla}\Vert_{L^2}\,\Vert\delta^{\nabla}R^{\nabla}\Vert_{L^4}{}^2\nonumber\\
&=\left(\frac{\sqrt{c_1}}{18}-2\Vert R^{\nabla}\Vert_{L^2}
\right)\,\Vert\delta^{\nabla}R^{\nabla}\Vert_{L^4}{}^2\nonumber\\
&\qquad+\left(
k-\frac19\left(\frac{c_1}{{\rm Vol}(M,g)}\right)^{1/2}
\right)\,\Vert\delta^{\nabla}R^{\nabla}\Vert_{L^2}{}^2.
\end{align}
Since 
$\Vert\delta^{\nabla}R^{\nabla}\Vert_{L^2}{}^2\geq 0$ in (10.14), 
we also have 
\begin{equation}
{\rm the\,\, RHS\,\,of\,\,(10.14)}
\geq 
k\Vert \delta^{\nabla}R^{\nabla}\Vert_{L^2}{}^2-2\Vert R^{\nabla}\Vert_{L^2}\,\Vert\delta^{\nabla}R^{\nabla}\Vert_{L^4}{}^2.
\end{equation}
\par
\underline{Case 1}:\,\, $\Vert\delta^{\nabla}R^{\nabla}\Vert_{L^2}{}^2
\geq \frac{{\rm Vol}(M,g)^{1/2}}{2}
\Vert \delta^{\nabla}R^{\nabla}\Vert_{L^4}{}^2$.
\qquad
In this case, if 
$\Vert\delta^{\nabla}R^{\nabla}\Vert_{L^4}>0$, then 
\begin{align}
{\rm the\,\, RHS \,\,of \,\,(10.17)} &> k\Vert\delta^{\nabla}R^{\nabla}\Vert_{L^2}{}^2
-\frac{k}{2}\,{\rm Vol}(M,g)^{1/2}\Vert\delta^{\nabla}R^{\nabla}\Vert_{L^4}{}^2\nonumber\\ 
&\qquad\qquad\qquad
\left({\rm by\,\,}2\Vert R^{\nabla}\Vert_{L^2}<\frac{k}{2}{\rm Vol}(M,g)^{1/2}\right)\nonumber\\
&=k\left(\Vert \delta^{\nabla}R^{\nabla}\Vert_{L^2}{}^2-\frac{{\rm Vol}(M,g)^{1/2}}{2}\Vert\delta^{\nabla}R^{\nabla}\Vert_{L^4}{}^2\right)\nonumber\\
&\geq 0\nonumber
\end{align} 
which is a contradication. We have $\Vert\delta^{\nabla}R^{\nabla}\Vert_{L^4}=0$, i.e., $\delta^{\nabla}R^{\nabla}=0$. 
\par
\underline{Case 2.}:
$\Vert\delta^{\nabla}R^{\nabla}\Vert_{L^2}{}^2
\leq \frac{{\rm Vol}(M,g)^{1/2}}{2}
\Vert \delta^{\nabla}R^{\nabla}\Vert_{L^4}{}^2$. 
\qquad
In this case, if 
$\Vert \delta^{\nabla}R^{\nabla}\Vert_{L^2}>0$, 
then 
\begin{align}
{\rm the\,\, RHS\,\,of\,\,(10.16)}
&=
\left(\frac{\sqrt{c_1}}{18}-2\Vert R^{\nabla}\Vert_{L^2}
\right)\,\Vert\delta^{\nabla}R^{\nabla}\Vert_{L^4}{}^2\nonumber\\
&\qquad+\left(
k-\frac19\left(\frac{c_1}{{\rm Vol}(M,g)}\right)^{1/2}
\right)\,\Vert\delta^{\nabla}R^{\nabla}\Vert_{L^2}{}^2\nonumber\\
&\geq 
\left(
\frac{\sqrt{c_1}}{18}-2\Vert R^{\nabla}\Vert_{L^2}
\right)\,2\,{\rm Vol}(M,g)^{-1/2}\Vert\delta^{\nabla}R^{\nabla}\Vert_{L^2}{}^2
\nonumber\\
&\qquad+\left(
k-\frac19\left(\frac{c_1}{{\rm Vol}(M,g)}\right)^{1/2}
\right)\,\Vert\delta^{\nabla}R^{\nabla}\Vert_{L^2}{}^2\nonumber\\
&\qquad\qquad\qquad\left(
{\rm by}\,\,\frac{\sqrt{c_1}}{18}-2\Vert R^{\nabla}\Vert_{L^2}\geq 0
\right)\nonumber\\
&=
\bigg\{
\frac{\sqrt{c_1}}{9}\,{\rm Vol}(M,g)^{-1/2}
-2\Vert R^{\nabla}\Vert_{L^2}\cdot 2{\rm Vol}(M,g)^{-1/2}
\nonumber\\
&\qquad+k-\frac19
\left(\frac{c_1}{{\rm Vol}(M,g)}
\right)^{1/2}
\bigg\}\,\Vert\delta^{\nabla}R^{\nabla}\Vert_{L^2}{}^2
\nonumber\\
&=\left(
k-2\Vert R^{\nabla}\Vert_{L^2}\,\cdot\,2\,{\rm Vol}(M,g)^{-1/2}
\right)\,\Vert\delta^{\nabla}R^{\nabla}\Vert_{L^2}{}^2\nonumber\\
&>0,\nonumber
\end{align}
which is also a contradiction. Thus, we have 
$\Vert\delta^{\nabla}R^{\nabla}\Vert_{L^2}=0$, i.e., 
$\delta^{\nabla}R^{\nabla}=0$. 
\qed
\vskip1.8cm\par

\end{document}